\documentclass{amsart}
\usepackage{amssymb,amscd}
\usepackage{enumerate}
\begin{document}

\setcounter{tocdepth}{1}
%
%
\setlength{\marginparwidth}{1.12in}

%
%
\newcommand{\wt}[1]{\widetilde{#1}}
\newcommand{\Prod}{\prod}
\newcommand{\Cal}{\mathcal}
\newcommand\Ran{\operatorname{Im}}
\newcommand\abs{\operatorname{abs}}
\newcommand\rel{\operatorname{rel}}
\newcommand\inv{\operatorname{inv}}
\newcommand\topo{\operatorname{top}}
\newcommand\opp{\operatorname{op}}
\newcommand\mfk{\mathfrak}
\newcommand\coker{\operatorname{coker}}
\newcommand\hotimes{\hat \otimes}
\newcommand\ind{\operatorname{ind}}
\newcommand\End{\operatorname{End}}
\newcommand\per{\operatorname{per}}
\newcommand\pa{\partial}
\newcommand\sign{\operatorname{sign}}
\newcommand\supp{\operatorname{supp}}
\newcommand\cy{\mathcal{C}^\infty}
\newcommand\CI{\mathcal{C}^\infty}
\newcommand\CO{\mathcal{C}_0}
\newcommand\lra{\longrightarrow}
\newcommand\vlra{-\!\!\!-\!\!\!-\!\!\!\!\longrightarrow}
\newcommand\bS{{}^b\kern-1pt S}
\newcommand\bT{{}^b\kern-1pt T}
\newcommand\Hom{\operatorname{Hom}}
\newcommand\LD{{\mathcal D}}

\newcommand{\cA}{\mathcal{A}}
\newcommand{\cB}{\mathcal{B}}
\newcommand{\cC}{\mathcal{C}}
\newcommand{\cD}{\mathcal{D}}
\newcommand{\cH}{\mathcal{H}}
\newcommand{\cJ}{\mathcal{J}}
\newcommand{\cJi}{\cJ_{\infty}}
\newcommand{\cJs}{\cJ^{*}}
\newcommand{\cun}{\cC^{\infty}}
\newcommand{\cunb}{\cC^{\infty}_{b}}
\newcommand{\cunc}{\cun_{c}}
\newcommand{\cuncd}{\dot{\cC}^{\infty}_{c}}
\newcommand{\cund}{\dot{\cC}^{\infty}}
\newcommand{\fA}{\mathfrak{A}}
\newcommand{\fAi}{\fA_{-\infty}}
\newcommand{\id}{{\rm id}}
\newcommand{\lzmn}{L^{2}(M_{0})}
\newcommand{\norm}[2]{\|#1\|_{#2}}
\newcommand{\nzn}{\NN_{0}}
\newcommand{\ons}{\setminus\{0\}}
\newcommand{\otl}{\omega_{T}^{\ell}}
\newcommand{\otlr}{\omega_{T}^{\ell,r}}
\newcommand{\otr}{\omega_{T}^{r}}
\newcommand{\RC}{{\rm RC}}
\newcommand{\rpq}{\overline{\RR}_{+}}
\newcommand{\sh}[1]{#1^{\sharp}}


\newcommand\alg[1]{\mathfrak{A}(#1)}
\newcommand\qalg[1]{\mathfrak{B}(#1)}
\newcommand\ralg[1]{\mathfrak{A}_r(#1)}
\newcommand\rqalg[2]{\mathfrak{B}_r(#1)}
\newcommand\ideal[1]{C^*(#1)}
\newcommand\rideal[1]{C^*_r(#1)}
\newcommand\qideal[2]{\mathfrak{R}_{#1}(#2)}
\newcommand\In{\operatorname{In}}

\newcommand\TR{\operatorname{T}}
\newcommand\ha{\frac12}
\newcommand\cal{\mathcal}
\newcommand\END{\operatorname{END}}
\newcommand\ENDG{\END_{\GR}(E)}
\newcommand\GL{\operatorname{GL}}
\newcommand\SU{\operatorname{SU}}
\newcommand\Cl{\operatorname{Cl}}
\newcommand\CC{\mathbb C}
\newcommand\NN{\mathbb N}
\newcommand\RR{\mathbb R}
\newcommand\ZZ{\mathbb Z}
\newcommand\ci{${\mathcal C}^{\infty}$}
\newcommand\CIc{{\mathcal C}^{\infty}_{\text{c}}}
\newcommand\hden{{\Omega^{\lambda}_d}}
\newcommand\VD{{\mathcal D}}
\newcommand\mhden{{\Omega^{-1/2}_d}}
\newcommand\ehden{r^*(E)\otimes {\Omega^{\lambda}_d}}


\newcommand{\Cat}{\mathcal C}
\newcommand{\Gr}[1]{{\mathcal G}^{(#1)}}
\newcommand{\GR}{\mathcal G}
\newcommand{\LGR}{\mathcal L}
\newcommand{\BB}{\mathbb{B}}
\newcommand{\GG}{\mathcal G}
\newcommand{\OA}{\mathcal O}
\newcommand{\PS}[1]{\Psi^{#1}(\GR;E)}
\newcommand{\tPS}[1]{\Psi^{#1}(\GR)}
\newcommand{\tPSI}[1]{\Psi^{#1}(M,\GR)}
\newcommand{\PSI}[1]{\Psi^{#1}(M,\GR;E)}
\newcommand{\AL}{{\mathcal A}(\GR)}
\newcommand{\FAM}{P=(P_x,x \in \Gr0)}
\newcommand\symb[2]{{\mathcal S}^{#1}(#2)}
\newcommand{\loc}{\operatorname{loc}}
\newcommand{\cl}{\operatorname{cl}}
\newcommand{\A}{s}
\newcommand{\prop}{\operatorname{prop}}
\newcommand{\comp}{\operatorname{comp}}
\newcommand{\adb}{\operatorname{adb}}
\newcommand{\dist}{\operatorname{dist}}

\newcommand{\alp }{r }
\newcommand{\bet }{d }
\newcommand{\gm }{\Gamma }
\newcommand{\lon }{\longrightarrow }
\newcommand{\be }{\begin{eqnarray*}}
\newcommand{\ee }{\end{eqnarray*}}
\newcommand{\GGR}{{\GR}}
\newcommand{\cald}{{\cal D}}
\newcommand{\calx}{{\cal X}}
\newcommand{\II}{\CIc(S^*(\GR), \End(E) \otimes {\mathcal P}_m)}
\newcommand{\IIY}{\CIc(S^*(\GR\vert_Y),
\End(E\vert_Y) \otimes {\mathcal P}_m)}
\newcommand{\mI}{\mathfrak I}
\newcommand\mF{\mathcal F}
\newcommand\mE{\mathcal E}
\newcommand\mD{\mathcal D}
\newcommand\mH{\mathcal H}
\newcommand\mL{\mathcal L}
\newcommand\Dir{{\not \!\!D}} 
\newcommand{\Cliff}{{\rm Cliff}}

\def\nin{\noindent}
\def\eg{e.g.\ }
\def\pt#1#2{\frac{\partial #1}{\partial #2}}

\newcommand{\frakg}{{\mathfrak g}}

\let\Tilde=\widetilde
\let\Bar=\overline
\let\Vec=\overrightarrow
\let\ceV=\overleftarrow
\def\vlra{\hbox{$\,-\!\!\!-\!\!\!-\!\!\!-\!\!\!-\!\!\!
-\!\!\!-\!\!\!-\!\!\!-\!\!\!-\!\!\!\longrightarrow\,$}}

\def\vleq{\hbox{$\,=\!\!\!=\!\!\!=\!\!\!=\!\!\!=\!\!\!
=\!\!\!=\!\!\!=\!\!\!=\!\!\!=\!\!\!=\!\!\!=\!\!\!=\!\!\!=\,$}}

\def\lrah{\hbox{$\,-\!\!\!-\!\!\!
-\!\!\!-\!\!\!-\!\!\!-\!\!\!-\!\!\!\longrightarrow\,$}}

\def\surj{-\!\!\!-\!\!\!-\!\!\!\gg}

\def\inj{>\!\!\!-\!\!\!-\!\!\!-\!\!\!>}

\let\<\langle
\let\>\rangle

\let\rho\varrho
\let\Ga\Gamma
\def\Gavert{{\Gamma_{\textrm{vert}}}}
\def\eref#1{(\ref{#1})}


 %
%
\newcommand\Mand{\text{ and }}
\newcommand\Mandset{\text{ and set}}
\newcommand\Mas{\text{ as }}
\newcommand\Mat{\text{ at }}
\newcommand\Mfor{\text{ for }}
\newcommand\Mif{\text{ if }}
\newcommand\Min{\text{ in }}
\newcommand\Mon{\text{ on }}
\newcommand\Motwi{\text{ otherwise }}
\newcommand\Mwith{\text{ with }}
\newcommand\Mwhere{\text{ where }}
\newcommand\ie{{i.\thinspace e.{}, }}
\newcommand\VV{\mathcal V}
\newcommand{\Diff}[1]{{\rm Diff}(#1)}
%
%
\newtheorem{theorem}{Theorem}[section]
\newtheorem{proposition}[theorem]{Proposition}
\newtheorem{corollary}[theorem]{Corollary}
\newtheorem{conjecture}[theorem]{Conjecture}
\newtheorem{lemma}[theorem]{Lemma}
\newtheorem{notation}[theorem]{Notations}
\theoremstyle{definition}
\newtheorem{definition}[theorem]{Definition}
\theoremstyle{remark}
\newtheorem{remark}[theorem]{Remark}
\newtheorem{example}[theorem]{Example}

%
%
\def\ie{{i.\thinspace e.{}, }}
\def\ppso{{\pi_{\mbox{\scriptstyle PSO}}}}
\def\PSO{\operatorname{P_{\scriptstyle SO}}}
\def\SO{\operatorname{SO}}
\def\so{{so}}
\def\Spin{\operatorname{Spin}}
\def\Ad{\operatorname{Ad}}
\newcommand{\mnote}[1]{\marginpar{\tiny\em #1}}

\author[B. Ammann]{Bernd Ammann} \address{Universit\"at Hamburg,
Fachbereich 11 -- Mathematik, Bundesstra\ss{}e 55, D - 20146 Hamburg, Germany}
\email{ammann@math.uni-hamburg.de}

\author[R. Lauter]{Robert Lauter} \address{Universit\"at Mainz,
Fachbereich 17 -- Mathematik, D - 55099 Mainz, Germany}
\email{lauter@mathematik.uni-mainz.de}

\author[V. Nistor]{Victor Nistor} \address{Pennsylvania State
       University, Math. Dept., University Park, PA 16802}
       \email{nistor@math.psu.edu}

\thanks{Nistor was partially supported by NSF Grants DMS-9971951,
DMS 02-00808, and ``collaborative research'' grant.
Ammann was partially supported
by The European Contract Human Potential Programme,
Research Training Networks HPRN-CT-2000-00101 and HPRN-CT-1999-00118.
}

\dedicatory\date
\begin{abstract}\
A manifold with a ``Lie structure at infinity'' is a non-compact
manifold $M_0$ whose geometry is described by a compactification
to a manifold with corners $M$ and a Lie algebra $\VV$ of vector
fields on $M$ subject to constraints only on $M \smallsetminus
M_0$.
This definition recovers several classes of non-compact manifolds
that were studied before: manifolds with cylindrical ends,
manifolds that are Euclidean at infinity, conformally compact
manifolds, and others. It hence provides a unified setting for the
study of these classes of manifolds and of their geometric
differential operators.
The Lie structure at infinity on $M_0$ determines a complete
metric on $M_0$ up to bi-Lipschitz equivalence.  This leads to the
natural problem of understanding the Riemannian geometry of these
manifolds, which is the main question addressed in this paper. We
prove, for example, that on a manifold with a Lie structure at
infinity the curvature tensor and its covariant derivatives are
bounded, by extending the Levi-Civita connection to an
$A^*$-valued connection where the bundle $A$ is uniquely
determined by the Lie algebra $\VV$.  We study a generalization of
the geodesic spray and give conditions for these manifolds to have
positive injectivity radius. We also prove that the geometric
operators are generated by the given Lie algebra of vector fields.

An important motivation for our study is to prepare the ground for
the investigation of the analysis of geometric operators on
manifolds with a Lie structure at infinity. The simplest examples
of manifolds with a Lie structure at infinity are the manifolds
with cylindrical ends. For these manifolds the corresponding
analysis is that of totally characteristic operators on a compact
manifold with boundary equipped with a ``$b$-metric.''
\end{abstract}

\title[Riemannian manifolds with a Lie structure at infinity]{On the
geometry of Riemannian manifolds with a Lie structure at infinity}

\maketitle \tableofcontents

\def\fr{\frac}
\def\ub{\underbar}
\def\O{\Cal O}
\def\F{\Cal F}
\def\differ{\text{differentiable} }
\def\tPSeudo{\text{pseudodifferential} }
\def\supp{\text{supp} }
\def\inn{{\mathcal R}}
\def\prop{\text{prop}}
\def\frag{\frak{G}}
\def\simd{\tilde{d}}
\def\simf{\tilde{\F}}
\def\simo{\tilde{\O}}
\def\simr{\tilde{r}}
\def\simp{\tilde{p}}
\def\simmu{\tilde{\mu}}

\def\thickpull#1{{#1}^{\#}}
\long\def\komment#1{}


\def\al{{\alpha}}
\def\be{{\beta}}
\def\de{{\delta}}
\def\De{{\Delta}}
\def\om{{\omega}}
\def\Om{{\Omega}}
\def\la{{\lambda}}
\def\ka{{\kappa}}
\def\io{{\iota}}
\def\si{{\sigma}}
\def\Si{{\Sigma}}
\def\ga{{\gamma}}
\def\ep{{\varepsilon}}
\def\Ga{{\Gamma}}
\def\Tau{\mathcal{T}}
\def\th{{\vartheta}}
\def\Th{{\Theta}}
\def\ph{{\varphi}}
\def\phi{{\varphi}}
\def\Ph{{\Phi}}
\def\rh{{\rho}}
\def\ze{{\zeta}}

\def\na{{\nabla}}
\def\pa{{\partial}}
\def\el{{\ell}}


\def\vol{{\mathop{{\rm vol}}}}
\def\area{{\mathop{{\rm area}}}}

\section*{Introduction}

Geometric differential operators on complete, non-compact Riemannian
manifold were extensively studied due to their applications to
physics, geometry, number theory, and numerical analysis.  Still,
their properties are not as well understood as those of differential
operators on compact manifolds, one of the main reason being that
differential operators on non-compact manifolds do not enjoy some of
the most useful properties enjoyed by their counterparts on compact
manifolds.

For example, elliptic operators on non-compact manifolds are not
Fredholm in general.  (We use the term ``elliptic'' in the sense
that the principal symbol is invertible outside the zero section.)
Also, one does not have a completely satisfactory
pseudodifferential calculus on an arbitrary complete, non-compact
Riemannian manifold, which might allow us to decide whether a
given geometric differential operator is bounded, Fredholm, or
compact
(see however \cite{aln2} and the references within). \

However, if one restricts oneself to certain classes of complete,
non-compact Riemannian manifolds, one has a chance to obtain more
precise results on the analysis of the geometric differential
operators on those spaces. This paper is the first in a series of
papers devoted to the study of such a class of Riemannian
manifolds, the class of Riemannian manifolds with a ``Lie
structure at infinity'' (see Definition \ref{def.unif.str}). We
stress here that few results on the geometry of these manifolds
have a parallel in the literature, although there is a fair number
of papers devoted to the analysis on {\em particular} classes of
such manifolds \cite{ConnesFOL, Cordes1, CordesDoong, CordesOwen1,
CordesOwen2, Kondratiev, Kordyukov, MelroseCongress,
MelroseScattering, Parenti, Roe, Schrohe, schulze01, Shubin,
TaylorG, WeinsteinAMS}. The philosophy of Cordes' comparison
algebras \cite{Cordes2}, Kondratiev's approach to analysis on
singular spaces \cite{Kondratiev}, Parenti's work on manifolds
that are Euclidean at infinity \cite{Parenti}, and Melrose's
approach to pseudodifferential analysis on singular spaces
\cite{MelroseCongress} have played an important role in the
development of this subject.

A manifold $M_0$ with a Lie structure at infinity has, by definition,
a natural compactification to a manifold with corners $M = M_0 \cup
\pa M$ such that the tangent bundle $TM_0 \to M_0$ extends to a vector
bundle $A \to M$ with some additional structure. We assume, for
example, that the Lie bracket of vector fields on $M_0$ defines, by
restriction, a Lie algebra structure on the space of sections of $A$
such that the space $\VV := \Gamma(A)$ of sections of $A$ identifies
with a Lie subalgebra of the Lie algebra of all vector field on
$M_0$. The pair $(M,\VV)$ then defines a Lie structure at infinity on
$M_0$.  A simple, non-trivial class of manifolds with a Lie structure
at infinity is that of manifolds with cylindrical ends. Let $M_0$ be a
manifold with cylindrical ends.  In this case, the compactification
$M$ is a manifold with boundary, $\VV$ consists of all vector fields
tangent to the boundary of $M$.  This example plays a prominent role
in the analysis of boundary value problems on manifolds with conical
points \cite{Kondratiev, MazyaBook, MitreaNistor, meaps, ss1, ss2}.
See the above references for earlier results.

Let $(M,\VV)$ be a Lie structure at infinity on $M_0$, $\VV =
\Gamma(A)$.  The choice of a fiberwise scalar product on $A$ gives
rise to a fiberwise scalar product $g$ on $TM_0$, i.e.\ a Riemannian
metric on $M_0$.  Since $M$ is compact, any two such metrics $g_1$ and
$g_2$ are equivalent, in the sense that there exists a positive
constant $C > 0$ such that $C^{-1}g_1 \ge g_2 \ge Cg_1$. One can thus
expect that the properties of the Riemannian manifold $(M_0,g)$
obtained by the above procedure depend only on the Lie structure at
infinity on $M_0$ and not on the particular choice of a metric on~$A$.
However, as shown in the following example, a metric on $M_0$ does not
determine a Lie structure at infinity on $M_0$.

\begin{example}
Let us compactify $\RR$ by including $+\infty$ and $-\infty$:
$$
        \overline{\RR} :=\RR \cup \{+\infty,-\infty\}.
$$
We define $\phi:[-1,+1]\to\overline{\RR}$, $\phi(t)=
\log(t+1)-\log(1-t)$, $\phi(\pm 1)=\pm \infty$. The pullback of the
differentiable structure on $[-1,1]$ defines a differentiable
structure on $\overline{\RR}$.  On $\overline{\RR}$ we consider the
Lie algebra of vector fields that vanish at $\pm \infty$.  The product
of these compactifications of $\RR$ defines a Lie structure on $M_0 :=
\RR^n$, in which the compactification $M$ is diffeomorphic to the
manifold with corners $[-1,1]^n$ and the sections of $A$ are all
vector fields tangent to all hyperfaces. (The resulting Lie structure
at infinity is that of the $b$-calculus (see Example \ref{bcal})).
Alternatively, one can consider the radial compactification of
$\RR^n$. The resulting Lie structure at infinity is described in
Example \ref{scat}, which is closely related to the so-called
scattering calculus \cite{MelroseScattering, Parenti}.

We thus see that $\RR^n$ fits into our framework and is in fact a
manifold with a Lie structure at infinity for several distinct
compactifications $M$.
\end{example}

Thus, although our motivation for studying manifolds with a Lie
structure at infinity comes from analysis, this class of manifolds
leads to some interesting questions about their geometry, and this
paper (the first one in a series of papers on this subject) is devoted
mainly to the issues and constructions that have a strong
Riemannian geometric flavor. It is important to mention here that
only very few results on the geometry of particular classes of
Riemannian manifolds with a Lie structure at infinity were proved
before, except some special examples (e.g.\ compact manifolds and
manifolds with cylindrical ends).  For example, we prove that
$M_0$ is complete and has bounded curvature, in the sense that the
Riemannian curvature $R$ and all its covariant derivatives
$\nabla^kR$, with respect to the Levi-Civita connection, are
bounded. Also, under some mild assumptions on $(M,\VV)$, we prove
that $(M_0,g)$ has positive injectivity radius, and hence $M_0$
has bounded geometry.  This is very convenient for the analysis on
these manifolds. The main technique is based on generalizing the
Levi-Civita connection to an ``$A^*$-valued connection'' on $A$.
(An {\em $A^*$-valued connection on a bundle $E \to M$} is a
differential operator $\nabla : E \to E \otimes A^*$ that
satisfies all the usual properties of a connection, but with $A$
replacing the tangent bundle, see Definition \ref{def.AvalC}. This
concept was first introduced in a slightly different form in
\cite{EvensLuW} by Evens, Lu, and Weinstein.
The right approach to the geometry of manifolds with a Lie
structure at infinity requires us to replace the tangent bundle by
$A$. This was noticed before in particular examples, see for
instance \cite{ConnesFOL, LN1, MelroseScattering,
MelroseCongress, Monthubert1}.

The Lie structure at infinity on $M_0$ allows us to define a
canonical algebra of differential operators on $M_0$, denoted
$\Diff{\VV}$, as the algebra of differential operators generated
by the vector fields in $\VV = \Gamma(A)$ and multiplication by
functions in $\CI(M)$. If $E_0,E_1 \to M$ are vector bundles on
$M$, then one can similarly define the spaces $\Diff{\VV;E_0,E_1}$
(algebras if $E_0 = E_1$) of differential operators generated by
$\VV$ and acting on sections of $E_0$ with values sections of
$E_1$. All geometric operators on $M_0$ (de Rham, Laplace, Dirac)
will belong to one of the spaces $\Diff{\VV;E_0,E_1}$, for
suitable bundles $E_0$ and $E_1$. The proof of this result depends
on our extension of the Levi-Civita connection to an $A^*$-valued
connection.

Many questions in the analysis on non-compact manifolds or on the
asymptotics of various families of operators can be expressed in
terms of $\Diff{\VV}$. We refer to \cite{CordesOwen1,
CrainicFernandez, defr, Hormander1, LMN, LN1, Mazzeo,
MelroseScattering, MelroseCongress, MelroseNistor} for just a few
of the many possible examples in the literature.  Indeed, let
$\Delta = d^*d \in \Diff{\VV}$ be the scalar Laplace operator on
$M_0$.  Then $\Delta$ is essentially self-adjoint on
$\cunc(M_{0})$ by old results of Gaffney \cite{gaff} and Roelcke
\cite{roel} from 1951 and 1960.  Assume that $M_0$ has positive
injectivity radius, then $P(1 + \Delta)^{-m/2}$ and $(1 +
\Delta)^{-m/2}P$ are bounded operators on $L^2(M_0)$, for any
differential operator $P \in \Diff{\VV}$ of order at most $m$.
Cordes \cite{Cordes1, Cordes2} defined the comparison algebra
$\mfk A(M,\VV)$ as the norm closed algebra generated by the
operators $P(1 + \Delta)^{-m/2}$ and $(1 + \Delta)^{-m/2}P$, with
$P$ a differential operator $P \in \Diff{\VV}$ of order at most
$m$. The comparison algebra is useful because it leads to criteria
for differential and pseudodifferential operators to be compact or
Fredholm between suitable Sobolev spaces \cite{LMN, ALNV1, aln2}.

We expect manifolds with a Lie structure at infinity and especially
the analytic tools (pseudodifferential and asymptotic analysis)
that we have established in \cite{ALNV1,aln2} to play an important
role for solving some problems in geometric analysis simultaneously
for a large class of manifolds. Indeed, in special cases of manifolds with 
a Lie structure at infinity
the solutions to quite a few
interesting problems in geometric analysis rely heavily on those methods. 
For instance, consider 
asymptotically
Euclidean manifolds, a special case of Example~\ref{scat}.  In
general relativity one is interested in finding solutions to the
Einstein equations whose spatial part is asymptotically Euclidean.
Integration of the first nontrivial coefficient in the asymptotic
development of the metric at infinity yields the so-called
``mass'' of the solution \cite{bartnik}. The positive mass theorem
states that any non-flat asymptotically Euclidean Riemannian
manifold with non-negative scalar curvature has positive mass. An
elegant proof of the positive mass theorem by Witten
\cite{parker.taubes} uses Sobolev embeddings on such manifolds.
The positive mass theorem provides the final step in the proof of
the Yamabe conjecture on compact manifolds \cite{schoen}: Any
conformal class on a compact manifold $M$ admits a metric with
constant scalar curvature.  In order to prove the conjecture in
the locally conformally flat case, one replaces the metric $g$ on
$M$ by a scalar-flat metric $u\cdot g$ on $M\setminus\{p\}$ where
$u$ is a function $u(x)\to \infty$ for $x\to p$, and a
neighborhood of $p$ provides the asymptotically euclidean end, and
one applies the positive mass conjecture to this. On most
non-compact manifolds, the Yamabe problem is still unsolved.
However, special cases have been solved, e.g.\ on manifolds with
cylindrical ends \cite{akutagawa.botvinnik}.

Both the geometry and the analysis of asymptotically hyperbolic
manifolds have been the subject of articles in general relativity
and the analysis of 3-manifolds, see \cite{Anderson,acd}. One can
prove rigidity theorems \cite{andersson.dahl} for asymptotically
hyperbolic ends, or existence results for asymptotically
hyperbolic Einstein metrics \cite{Lee}. Similar rigidity problems
for asymptotically complex hyperbolic ends are subject in
\cite{boualem.herzlich, chrusciel.herzlich, herzlich}.

Or take
the construction of manifolds with with special
holonomy $SU(m)$, $Sp(m)$ and $G_2$ 
where the analysis of weighted
function spaces on manifolds which are quasi asymtotically locally
euclidean \cite{joycebook,joyce.ale,joyce.qale}
has been used.

In summary, our present program will
lead to a unified approach to  the analysis on various types of
manifolds with a ``good'' asymptotic behavior at infinity.

We now discuss the contents of each section.  In
Section~\ref{Sec.Struc.Lie} we introduce and study structural Lie
algebras of vector fields and the equivalent concept of boundary
tangential Lie algebroids.  A structural Lie algebra of vector
fields on a manifold with corners $M$ gives rise to a canonical
algebra $\Diff{\VV}$ of differential operators.  We include
numerous examples.

Then, in Section \ref{Sec.Unif.Str}, we specialize to the case that
the constraints are only on the boundary.  This special case is called
a ``manifold with a Lie structure at infinity.''  The Lie structure at
infinity defines a Riemannian metric on the interior of the
manifold. This metric is unique up to bi-Lipschitz equivalence.
Hence, the Lie structures at infinity is a tool for studying a large
class of open Riemannian manifolds.  We are interested in the analysis
on such open manifolds.

Section \ref{Sec.3} is devoted to the study of the geometry of
Riemannian manifolds with a Lie structure at infinity. We will prove
that these manifolds are complete and have bounded curvature (together
with all its covariant derivatives). This depends on an extension of
the Levi-Civita connection to an $A^*$-valued connection, the
appropriate notion of connection in this setting.  Then we investigate
the question of whether a Riemannian manifold with a Lie structure at
infinity has positive injectivity radius.

In Section \ref{Sec.Geom.Op} we introduce Dirac and generalized
Dirac operators and prove that they belong to $\Diff{\VV;W}$,
where $W$ is a Clifford module. The same property is shared by all
geometric operators (Laplace, de Rham, signature) on the open
manifold $M_0$.

We thank Sergiu Moroianu for several discussions on the subject.

\section{Structural Lie algebras and Lie algebroids\label{Sec.Struc.Lie}}

We introduce in this section the concept of structural Lie
algebras of vector fields, which is then used to define manifolds
with a Lie structure at infinity.

\subsection{Projective modules\label{subs.pm}}\
In this subsection, we recall some well-known facts
about projective modules over $\CI(M)$, where $M$ is a compact
manifold, possibly with corners.

Let $V$ be a $\CI(M)$ module with module structure $\CI(M) \times
V \ni (f, v) \mapsto fv \in V$. Let $x \in M$ and denote by $\mathfrak
p_x$ the set of functions on $M$ that vanish at $x \in M$. Then
$\mathfrak p_x V$ is a complex vector subspace of $V$ and
$V/\mathfrak p_x V$ is called the {\em geometric fiber of $V$ at
$x$}. In general, the geometric fibers of $V$ are complex vector
spaces of varying dimensions.

A subset $S \subset V$ will be called a {\em basis} of $V$ if
every element $v \in V$ can be written uniquely as $v = \sum_{s
\in S} f_s s$, with $f_s \in \CI(M)$,
$\#\{s\in S \,|\, f_s\not = 0\}<\infty$.
(In our applications, $S$ will always be a finite set, so we will not
have to worry about this last condition.)
A $\CI(M)$-module is called \emph{free (with basis $S$)} if it has a basis $S$.
Unlike the general case, the geometric fibers of a free module have constant
dimension, equal to the number of elements in the basis $S$. Note
however, that if $f: V \to W$ is a morphism of free modules, the
induced map between geometric fibers {\em may have non-constant rank}.
For example, it is possible that $f$ is injective,
but the induced map on the geometric fibers is not injective on
all fibers. An example is provided by $M = [0, 1]$, $V = W =
\CI([0, 1])$ and $f$ being given by the multiplication with the
coordinate function $x \in [0, 1]$. Then $f$ is injective, but the
induced map on the geometric fibers at $0$ is $0$.

A $\CI(M)$--module $V$ is called {\em finitely generated
projective} if, by definition, there exists another module $W$
such that $V \oplus W$ is free with a finite basis. We then have
the following fundamental theorem of Serre and Swan \cite{Karoubi}

\begin{theorem}[Serre-Swan]\label{thm.SerreSwan}\
If $V$ is a finitely generated, projective module over $\CI(M)$,
then the set $E := \cup_{x \in M}  (V/ \mathfrak p_x V) \times
\{x\}$, the disjoint union of all geometric fibers of $V$, can be
endowed with the structure of a finite-dimensional, smooth vector
bundle $E \to M$ such that $V \simeq \Gamma(M; E)$. The converse
is also true:\ $\Gamma(M; E)$ is a finitely generated, projective
$\CI(M)$-module for any finite-dimensional, smooth vector bundle
$E \to M$.
\end{theorem}

Suppose now that $V$ is a $\CI(M)$-module and that $M$ is
connected. Then $V$ is a finitely generated, projective
$\CI(M)$-module if, and only if, there exists $k \in \ZZ_+$
satisfying the following condition:
\begin{quotation}
For any $x \in M$, there exist $\phi \in \CI(M)$, $\phi(x) = 1$,
and $k$-elements $v_1, \ldots, v_k \in V$ with the property that
for any $w \in V$ we can find $f_1, f_2, \ldots, f_k \in \CI(M)$
such that
\begin{equation}\label{eq.local.basis}
    \phi(f_1 v_1 + f_2 v_2 + \ldots + f_k v_k - w) = 0
    \text{ in } V
\end{equation}
and, moreover, the germs of $f_1, \ldots, f_k$ at $x$ are uniquely
determined.
\end{quotation}
A module $V$ satisfying Condition \eqref{eq.local.basis} above is
called {\em locally free of rank $k$}, and what we are saying here
is that ``locally free of rank $k$, for some $k$,'' is equivalent
to ``finitely-generated, projective.'' It is crucial here that the
number of elements $k$ is the same for any $x \in M$. In case $M$
is not connected, the number $k$ needs only be constant on the
connected components of $M$.

\begin{remark}
The introduction of projective modules over $C^\infty(M)$ in
Partial Differential Operators on non-compact manifolds was
pioneered by Melrose \cite{MelroseActa} in the early 1980s.
\end{remark}

\subsection{Manifolds with corners and structural Lie algebras}\
\label{subsec.cornerdef} We now fix our terminology and recall the
definitions of the main concepts related to manifolds with
corners.

In the following, by a {\em manifold} we shall always understand a
$C^\infty$-manifold {\em possibly with corners}. In contrast, a
{\em smooth manifold} is a $C^\infty$-manifold {\em without
corners}.  By definition, for every point $p$ in a manifold with
corners $M$, there is a coordinate neighborhood $U_p$ of $p$ and
diffeomorphism $\phi_p$ to $[0,\infty)^k \times \RR^{n-k}$, with
$\phi_p(p) = 0$, such that the transition functions are smooth
(including on the boundary). The number $k$ here clearly depends
on $p$, and will be called the \emph{boundary depth of $p$}. Hence
points in the interior have boundary depth $0$, points on the
boundary of a manifold without corners have boundary depth $1$,
etc. Roughly speaking the boundary depth counts the number of
boundary faces $p$ is in.

Moreover, we assume that each hyperface $H$ of $M$ is
an embedded submanifold and has a defining function, that is,
there exists a smooth function $x_H \ge 0$ on $M$ such that
$$
        H = \{ x_H = 0 \} \text{ and } dx_H \not = 0 \text{
        on }H.
$$ This assumption is just a simplifying assumption. We can deal
with general manifolds with corners using the constructions from
\cite{Monthubert1}. Note that a priori we do not fix a particular
system of defining functions, but only use their existence
occasionally.

If $F\subset M$ is an arbitrary face of $M$ of codimension $k$, then
$F$ is an open component of the intersection of the hyperfaces
containing it.
Any set $x_1,\ldots,x_k$ of defining functions of the
hypersurfaces containing $F$ is called a
{\em set of  defining functions of $F$};
thus, $F$ is a connected component of $\{x_1=x_2 =\ldots =x_k=0\}$.
This statement obviously does not depend on the choice of the
defining functions $x_j$, $j=1,\ldots,k$.
We
shall denote by $\pa M$ the union of all non-trivial faces of $M$.
Usually, we shall denote by $M_0$ the interior of $M$, that is,
$M_0:= M \smallsetminus \pa M$.

A {\em submersion} $f : M \to N$, between two manifolds with corners
$M$ and $N$, is a differentiable map $f$ such that $df$ is surjective
at all points and $df(v)$ is an inward pointing tangent vector of $N$
if, and only if, $v$ is an inward pointing vector $M$.  It follows
then that the fibers $f^{-1}(y)$ of $f$ are smooth manifolds {\em
without} corners, and that $f$ preserves the boundary depth (i.e.\ the
number of boundary faces a point $p$ is in.)  If $x$ is a defining
function of some hyperface of $N$, then $x_1 = x \circ f$ is such that
$\{x_1 = 0\}$ is a union of hyperfaces of $M$ and $dx_1 \not = 0$ on
$\{x_1 = 0\}$.

\begin{example}
Let $A\to M$ be a smooth vector bundle. The {\em sphere bundle of
$A$}, denoted $S(A)$, is defined, as usual, as the set of
(positive) rays in the bundle $A$, that is,
$S(A)=(A\ons)/\RR_{+}$. If we fix a smooth metric on $A$, then
$S(A)$ identifies with the set of vectors of length one in $A$.
Moreover, $S(A)\to M$ turns out to be a submersion of manifolds
with corners.
\end{example}

A {\em submanifold with corners} $N$ of a manifold with corners
$M$ is a submanifold $N \subset M$ such that $N$ is a manifold
with corners, and each hyperface $H$ of $N$ is a connected
component of a set of the form $H' \cap N$, where $H'$ is a
hyperface of $M$ intersecting $N$ {\em transversally}.

The starting point of our analysis is a Lie algebra of vector
fields on a manifold with corners. For reasons that will be
clearer later, we prefer to keep this concept as general as
possible, even if for the analysis on non-compact manifolds, only
certain classes of Lie algebras of vector fields will be used.

\begin{definition}\label{def.Lie.a.v}\
A {\em structural Lie algebra of vector fields} on a manifold $M$
($M$ possibly with corners) is a subspace $\VV \subset \Gamma(TM)$
of the real vector space of vector fields on $M$ with the
following properties:

(i)\ $\VV$ is closed under Lie brackets;

(ii)\ $\VV$ is a finitely generated, projective $\CI(M)$-module;
and

(iii)\ the vector fields in $\VV$ are tangent to all faces in $M$.
\end{definition}

By (ii) we mean that $\VV$ is closed for multiplication with
functions in $\CI(M)$ and the induced $\CI(M)$-module structure
makes it a finitely generated projective $\CI(M)$-module.

Given a structural Lie algebra $\VV$ of vector fields on a
manifold with corners, we call the enveloping algebra $\Diff{\VV}$
of $\VV$ the algebra of {\em $\VV$-differential operators on $M$}.
Note that any $\VV$-differential operator $P \in\Diff{\VV}$ can be
realized as a polynomial in vector fields in $\VV$ with
coefficients in $\cun(M)$ acting on the space $\cun(M)$

Let us give some examples for structural Lie algebras of vector
fields. Some of these examples can also be found in
\cite{MelroseScattering}. We also give descriptions of the
structural vector fields in local coordinates, because this will
be helpful in the applications of the theory developed here.  All
of the following examples model the analysis on some non-compact
manifold, except for the last one, which models the analysis of
adiabatic families.

The following example is the simplest and most studied so far,
however, it is quite important for us because it models the
geometry of manifolds with cylindrical ends, and hence it is
easier to grasp.

\begin{example}\
\label{bcal}
Let $M$ be a manifold with corners and
\begin{equation}\label{eq.def.Vb}
        \VV_b = \{X \in \Gamma(TM) : X \text{ is tangent to all
        faces of } M\,\}.
\end{equation}
Then $\VV_b$ is a structural Lie algebra of vector fields, and any
structural Lie algebra of vector fields on $M$ is contained in
$\VV_b$, by condition (iii) of the above definition. A vector
field $X\in\VV_{b}$ is called a {\em $b$-vector field $X$}. Fix
$x_{1},\ldots,x_{k}$ and $y\in\RR^{n-k}$ local coordinates near a
point $p$ on a boundary face of codimension $k$, with $x_j$
defining functions of the hyperfaces through $p$. Then any
$b$-vector field $X$ is of the form
\begin{equation*}
        X=\sum_{j=1}^{k}a_{j}(x,y)x_{j}\partial_{x_{j}} +
        \sum_{j=1}^{n-k}b_{j}(x,y)\partial_{y_{j}}
\end{equation*}
on some neighborhood of $p$, with the coefficients $a_{j}$ and $b_{j}$
smooth everywhere (including the hyperfaces $x_{j}=0$), for all
$j$. This shows that the Lie algebra of $b$-vector fields is generated
in a neighborhood $U$ of $p$ by $x\partial_{x}$ and $\partial_{y}$ as
a $\CI(M)$-module. The differential operators in $\Diff{\VV_b}$ are
called {\em Fuchs type operators}, {\em totally characteristic}, or
simply, and perhaps more systematically {\em $b$-differential
operators}.  The structural Lie algebra $\VV_{b}$ and the analysis of
the corresponding differential and pseudodifferential operators are
treated in detail for instance in \cite{eichorn1, eichorn5,
Hormander3, lesch, meaps, MelroseScattering, schwil}.
\end{example}

\begin{example}
\label{scat} Let $M$ be a compact manifold with boundary and
$x:M\rightarrow\rpq$ a boundary defining function. Then the Lie
algebra $\VV_{sc}:=x\VV_{b}$ does not depend on the choice of $x$
and the vector fields in $\VV_{sc}$ are called {\em scattering
vector fields}; with respect to local coordinates $(x,y)$ near the
boundary, scattering vector fields are generated by
$x^{2}\partial_{x}$ and $x\partial_{y}$.  An analysis of the
scattering structure can be found in \cite{MelroseScattering}.
Since this structure models the analysis on asymptotically
Euclidean spaces, let us be a little bit more precise and recall
some basic definitions.  A Riemannian metric $g$ on the interior
of $M$ is called a {\em scattering metric} if, close to the
boundary $\pa M$, it is of the form
$g=\frac{dx^{2}}{x^{4}}+\frac{h}{x^{2}}$ where $h$ is a smooth,
symmetric $2$-tensor on $M$ which is non-degenerate when
restricted to the boundary. Then scattering vector fields are
exactly those smooth vector fields on $M$ that are of bounded
length with respect to $g$, and the corresponding Laplacian
$\Delta_{g}$ is an elliptic polynomial in scattering vector
fields.  As a special case of this setting note that the radial
compactification map
\begin{equation*}
        \RC:\RR^{N}\longrightarrow
        S^{N}_{+}:=\{\omega=(\omega_{0},\omega')\in S^{N}:
        \omega_{0}\geq0\}: z\longmapsto (1+|z|^{2})^{-1/2}(1,z)
\end{equation*}
identifies the Euclidean space $\RR^{N}$ with the interior of the
upper half-sphere $S^{N}_{+}$ such that the Euclidean metric lifts to
a scattering metric on $S^{N}_{+}$.
\end{example}

The following example is one of the examples that we are
interested to use in applications.

\begin{example}\label{edge}\
Let $M$ be a manifold with boundary $\pa M$, which is the total
space of a fibration $\pi : \pa M \to B$ of smooth manifolds. We
let
\begin{equation*}
        \VV_e = \{X \in \Gamma(TM): X \text{ is tangent to all
        fibers of } \pi\, \text{ at the boundary}\}
\end{equation*}
be the space of {\em edge vector fields}.  In order to show that
this is indeed a structural Lie algebra of vector fields, we have
to show that it is closed under  Lie brackets.  Let $i:\pa M\to M$
be the inclusion.  Assume that $X,Y\in \VV_e$. Because
\begin{equation*}
        [X,Y]|_{\pa M}=[X|_{\pa M},Y|_{\pa M}],
\end{equation*}
the commutator is again tangent to the fibers of $\pi$. If
$(x,y,z)$ are coordinates in a local product decomposition near
the boundary, where $x$ corresponds to the boundary defining
function, $y$ to a set of variables on the base $B$ lifted through
$\pi$, and $z$ is a set of variables in the fibers of $\pi$, then
edge vector fields are generated by $x\partial_{x}$,
$x\partial_{y}$, and $\partial_{z}$. Using this local coordinate
description is another way, to see immediately that the space of
edge vector fields is in fact a Lie algebra. More importantly, it
shows that it is a projective $\CI(M)$-module. The analysis of the
Lie algebra $\VV_{e}$ is partly carried out in \cite{Mazzeo} and
more recently in \cite{LN1}.
\end{example}

A special case of the edge structure is of particular importance
for the analysis on hyperbolic space, so it deserves its own name:

\begin{example}
\label{zero} Let $M$ be a compact manifold with boundary, and let
$\VV_{0}$ be the edge vector fields corresponding to the trivial
fibration $\pi = \id: \partial M\rightarrow\partial M$, \ie we
have
\begin{equation*}
        \VV_{0}=\{X\in\Gamma(TM):X|_{\partial M}=0\}
\end{equation*}
which explains the name {\em $0$-vector fields} for the elements
in $\VV_{0}$.  With respect to local coordinates $(x,y)$ near the
boundary, $0$-vector fields are generated by $x\partial_{x}$ and
$x\partial_{y}$.  Recall that a Riemannian manifold
$(M_{0},g_{0})$ is called {\em conformally compact} provided it is
isometric to the interior of a compact manifold $M$ with boundary
equipped with a metric $g=\varrho^{-2}h$ in the interior, where
$h$ is a smooth metric on $M$ and $\varrho:M\rightarrow\rpq$ a
boundary defining function. Note that $0$-vector fields are the
smooth vector fields on $M$ that are of bounded length with
respect to $g$; moreover, the Laplacian $\Delta_{g_{0}}$ is given
as an elliptic polynomial in $0$-vector fields.  A particular
example of conformally compact spaces is of course the hyperbolic
space with compactification given by the ball model.  Conformally
compact spaces arise naturally in questions related to the
Einstein equation \cite{Anderson, Lee, MP}, and the
``AdS/CFT-correspondence.''  An analysis of $0$-vector fields and
the associated $0$-differential and pseudodifferential operators
was carried out for instance in \cite{zfr, Mazzeo, rade98}.
Criteria for the Fredholmness of operators in $\Diff{\VV_0}$,
which is crucial in the approach to the study of Einstein's
equations on conformally compact manifolds used in the above
mentioned papers, were established for instance in \cite{zfr, LMN,
LN0, LN1, Mazzeo, mame87, MelroseScattering, rade98}.
\end{example}

The structural Lie algebra of vector fields in the next example is
a slight variation of the Lie algebra of edge vector fields,
however, it is worth pointing out that this slight variation leads
to a completely different analysis for the associated
(pseudo)differential operators.

\begin{example}
\label{de} Let $M$ be as in Example \ref{edge} and
$x:M\rightarrow\rpq$ be a boundary defining function. Then
$\VV_{de}:=x\VV_{e}$ is a structural Lie algebra of vector fields;
the corresponding structure is called the {\em double-edge
structure}. With respect to local product coordinates as in
Example \ref{edge}, double-edge vector fields are generated by
$x^{2}\partial_{x}$, $x^{2}\partial_{y}$, and $x\partial_{z}$.
The analysis of the double-edge structure, which is in fact much
simpler than the corresponding analysis of the edge structure, can
be found for instance in \cite{defr}.
\end{example}

The following example appears in the analysis of holomorphic
functions of several variables.

\begin{example}
\label{theta} Let $M$ be a smooth compact manifold with boundary
$\pa M$ and let $\Theta\in\cun(M,\Lambda^{1}T^{*}M)$ be a smooth
$1$-form such that $i^{*}\Theta\neq0$ where $i:\partial
M\hookrightarrow M$ is the inclusion of the boundary. Moreover,
let $x$ be a boundary defining function. Then
\begin{equation*}
        \VV_{\Theta}:=\{V\in\VV_{b}: V=0 \text{ at }
        \partial M \text{ and } \Theta(V)\in x^{2}\cun(M)\}
\end{equation*}
is a structural Lie algebra of vector fields that is called the
{\em $\Theta$-structure}. For a local description as well as for
an analysis of the $\Theta$-structure we refer to \cite{emm91}.
\end{example}

All the above examples of structural Lie algebras of vector fields
model the analysis on certain non-compact manifolds (giving rise
to algebras of differential operators that replace the algebra of
totally characteristic differential operators) on manifolds with
cylindrical ends. The following example, however, models the
analysis of a family of an adiabatic differential operators.

\begin{example}
\label{adia} Let $N$ be a closed manifold that is the total space
of a locally trivial fibration $\pi:N\rightarrow B$ of closed
manifolds, let $TN/B\rightarrow N$ be the vertical tangent bundle,
and let $M:=N\times[0,\infty)_{x}$. Then
\begin{equation*}
        \VV_{a}:=\{V\in\Gamma(TM): V(x)\in TN \text{ for all }
        x\in[0,\infty) \text{ and } V(0)\in \Gamma(TN/B)\}
\end{equation*}
is a structural Lie algebra of vector fields that is called the
{\em adiabatic algebra}. If $(y,z)$ are local coordinates on $N$,
where again the set of variables $y$ corresponds to variables on
the base $B$ lifted through $\pi$, and $z$ are variables in the
fibers, then adiabatic vector fields are generated by
$x\partial_{y}$ and $\partial_{z}$. The adiabatic structure has
been studied and used for instance in \cite{mordis} and
\cite{moroianu:01}
\end{example}

We shall sometimes refer to a structural Lie algebra of vector
fields simply as {\em Lie algebra of vector fields}, when no
confusion can arise. Because $\VV$ is a finitely generated,
projective $\CI(M)$-module, using the Serre-Swan theorem
\cite{Karoubi} (recalled above, see Theorem \ref{thm.SerreSwan})
we obtain that there exists a vector bundle
\begin{equation}
        A = A_{\VV} \to M, \quad \text{such that } \; \VV \cong
        \Gamma(A_{\VV}),
\end{equation}
{\em naturally} as $\CI(M)$-modules. We shall identify from now on
$\VV$ with $\Gamma(A_{\VV})$. The following proposition is
due to Melrose.

\begin{proposition}\label{prop.vbundle}\
If $\VV$ is a structural Lie algebra of vector fields, then there
exists a natural vector bundle map $\varrho : A_{\VV} \to TM$ such
that the induced map $\varrho_\Gamma : \Gamma(A_{\VV}) \to
\Gamma(TM)$ identifies with the inclusion map.
\end{proposition}

\begin{proof}\ Let $m \in M$. Then the fiber $A_{m}$ at $m$ of
$A = A_{\VV} \to M$ identifies with $\VV/\mathfrak p_m\VV$, where
$\mathfrak p_m$ is the ideal of $\CI(M)$ consisting of smooth
functions on $M$ that vanish at $m$.  Recall now that $\VV$ consists
of vector fields on $M$. Then the map $A_m \to T_m M$ sends the class
of $X \in \VV$ to the vector $X(m) \in T_m M$.
\end{proof}

\begin{remark}\label{remark.super.pol}
The condition in Definition~\ref{def.Lie.a.v} that $\VV$ has to be
projective is essential. As an example consider $M=[0,1]$ and let
\begin{multline*}
        \VV:=\Bigl\{f(t)\partial_t\, : \;
         f:[0,1]\to \RR \text{ smooth}\,,\;\;f(1)=0\,,
        \; \text{ and } \\ t^k(d^m f/dt^m)\to 0
        \text{ as } t\to 0\mbox{ for all }k,m\in \NN\cup \{0\}.\Bigr\}.
\end{multline*}

Then $\VV$ is a $C^\infty(M)$-module.  However, $\VV$ is not a
projective $C^\infty(M)$-module, as we can see by contradiction.
Assume $\VV$ were projective. Then there is a bundle $A$ over
$[0,1]$ with $\VV=\Gamma(A)$. Let $s$ be a trivialization of $A$,
\ie $s(t)=f(t)\partial_t$ with $f$ as above. Hence $\tilde f(t) =
(1/t) f(t)$ also decays sufficiently fast; however
\begin{equation*}
        \Gamma(A)\not\ni\frac{1}{t}s(t)=\tilde f(t)\partial_t\in \VV.
\end{equation*}
\end{remark}

It is convenient for the following discussion to recall the
definition of a Lie algebroid.  General facts about Lie algebroids
can be found in \cite{CrainicFernandez, Mackenzie1} (a few basic
facts are also summarized in \cite{NWX}).

\begin{definition}
\label{Lie.Algebroid} A  {\em Lie algebroid} $A$ over a manifold
$M$ is a vector bundle $A$ over~$M$, together with a Lie algebra
structure on the space $\Gamma(A)$ of smooth sections of $A$ and a
bundle map $\varrho: A \rightarrow TM$, extended to a map
$\varrho_\Gamma : \Gamma(A) \to \Gamma(TM)$ between sections of
these bundles, such that
\begin{enumerate}[(i)]
\item
$\varrho_\Gamma([X,Y])=[\varrho_{\Gamma}(X),\varrho_{\Gamma}(Y)]$,
\item $[X, fY] = f[X,Y] + (\varrho_{\Gamma}(X) f)Y$,
\noindent for any smooth sections $X$ and $Y$ of $A$ and any
smooth function $f$ on $M$.
\end{enumerate}
The map $\varrho_\Gamma$ is called the {\em anchor of} $A$. If, in
addition,
\begin{enumerate}[(i)]
\setcounter{enumi}{2}
\item all vector fields $\varrho_\Gamma(\Gamma(A))$ are tangential to the
faces,
\end{enumerate}
then the Lie algebroid $A \to M$ is called a \emph{boundary
tangential} Lie algebroid.
\end{definition}

We thus see that there exists an equivalence between the concept
of a structural Lie algebra of vector fields $\VV =
\Gamma(A_{\VV})$ and the concept of a boundary tangential Lie
algebroid $\varrho : A \to TM$ such that $\varrho_\Gamma :
\Gamma(A) \to \Gamma(TM)$ is injective and has range in $\VV_b$.
In order to shorten our notation, we will write $Xf$ instead of
$\varrho_\Gamma(X) f$ for the action of the sections of a Lie
algebroid on functions if the meaning is clear from the context.

\subsection{Constructing new Lie algebroids from old ones}
Let $f:N\to M$ be a submersion of manifolds with corners in the
above sense (which implies in particular that any fiber is a
smooth manifold).  Let $A=A_\VV$ be a boundary tangential Lie
algebroid over $M$.

\begin{definition}
The \emph{thick pull-back} $\thickpull{f} A$ is the vector bundle
over $N$ which at the point $p\in N$ is defined to be
\begin{equation*}
        \thickpull{f} A_p:= \left\{(v,w)\,|\,v\in A_{f(p)},\,
        w\in T_pN,\, f_*(w)=\rho(v)\right\}
\end{equation*}
equipped with the vector bundle structure induced by $f^*(A)
\oplus TN$.

Projection to the first component yields a surjective linear map
$\thickpull{f} A_p\to A_{f(p)}$, denoted in the following by
$f_*$, and projection onto the second component yields a linear
map $\thickpull{f} A_p\to T_pN$, denoted by $\thickpull{f} \rho$.
\end{definition}
We obtain the commuting diagram
\[
        \begin{array}{ccc} \thickpull {f} A_p &
        \stackrel{\thickpull{f} \rho}{ \lon } & T_pN \\[3mm]
        \downarrow f_* && \downarrow f_* \\[3mm] A_{f(p)} &
        \stackrel{\rho}{\lon } & T_{f(p)}M.  \end{array}
\]
For example $\thickpull{f}TM=TN$.

\begin{lemma}
The thick pull-back $f^{\#}A$ is a boundary tangential
Lie algebroid over $N$ with
anchor map given by $f^{\#}\varrho$.
\end{lemma}

\begin{proof}
Let $\Gavert{TN}$ denote the bundle of vertical sections $X$, \ie
$f _*X = 0$.
This bundle coincides by definition with the analogously
defined bundle of vertical sections of $A$.
The rows of the following commutative diagram are exact.
\[
        \begin{array}{ccccccccc} 0 & \lon & \Gavert(\thickpull{f}A) &
        \lon & \Ga(\thickpull{f}A) & \lon & \Ga(A)
        \otimes_{C^\infty(M)} C^\infty (N) & \lon & 0\\[3mm]
        &&\downarrow\cong && \downarrow && \downarrow \\[3mm] 0 & \lon
        & \Gavert(TN) & \lon & \Ga(TN) & \lon & \Ga(TM)
        \otimes_{C^\infty(M)} C^\infty (N) & \lon & 0\\ \end{array}
\]
The vertical arrows are inclusions.  The horizontal arrows of the
second row are Lie algebra homomorphisms.  The space $\Ga(A)$ is
by definition a Lie subalgebra of $\Ga(TM)$, thus
$\Ga(A)\otimes_{C^\infty(M)} C^\infty (N)$ is a Lie subalgebra of
$\Ga(TM) \otimes_{C^\infty(M)} C^\infty (N)$.  A standard diagram
chase then implies that $\Ga(\thickpull{f}A)$ is also a Lie
subalgebra of $\Ga(TN)$.

The fact that $A$ is projective [respectively, boundary
tangential] immediately implies that $\Ga(\thickpull{f}A)$ is also
projective [respectively, boundary tangential].
\end{proof}

Let $\mfk g$ and $\mfk h$ be two Lie algebras. Suppose that there
is given an action by derivations of $\mfk g$ on $\mfk h$:
\begin{equation}
        \phi : \mfk g \to \operatorname{Der}(\mfk h).
\end{equation}
Then we can define the semi-direct sum  $\mfk g \ltimes_\phi \mfk
h$ as follows. As a vector space, $\mfk g \ltimes_\phi \mfk h =
\mfk g \oplus \mfk h$, and the Lie bracket is given by
\begin{equation}
        [(X_1,Y_1), (X_2,Y_2)] := ([X_1,X_2], \phi_{X_1}(Y_2) -
        \phi_{X_2}(Y_1) + [Y_1,Y_2]),
\end{equation}
for any $X_1,X_2 \in \mfk g$ and $Y_1,Y_2 \in \mfk h$. We shall
usually omit the index $\phi$ denoting the action by derivations
in the notation for the semi-direct sum.

We want to use this construction to obtain new Lie algebroids from
old ones. Assume then that we are given two Lie algebroids $A, L
\to M$ over the same manifold and that $\Gamma(A)$ acts by
derivations on $\Gamma(L)$. Denote this action by $\phi$, as
above. We assume that the action of $\Gamma(A)$ on $\Gamma(L)$ is
compatible with the $\CI(M)$-module structure on $\Gamma(L)$, in
the sense that
\begin{equation*}
        \phi_X(fY) = X(f) Y + f \phi_X(Y),
\end{equation*}
for any $X \in \Gamma(A)$, $Y \in \Gamma(L)$, and $f \in \CI(M)$.
Assume, for simplicity, that the structural map (anchor)
$L \to TM$ is zero, then we can
endow $\Gamma(A \oplus L) = \Gamma(A) \oplus \Gamma(L)$ with the
semi-direct sum structure obtained from $\Gamma(A) \ltimes
\Gamma(L)$ such that $A \oplus L$ becomes a Lie algebroid denoted $A
\ltimes L$, and called the {\em semi-direct product of $L$ by $A$}
\cite{HigginsMackenzie}. Thus
\begin{equation}\label{eq.def.cr.p}
        \Gamma(A \ltimes L) = \Gamma(A) \ltimes \Gamma(L).
\end{equation}

In the language of Lie algebroids, the action of $\Gamma(A)$ on
$\Gamma(L)$ considered above is called a representation of $A$ on
$L$. In a similar way, the action of $\Gamma(A)$ on $\Gamma(L)$ {\em
by derivation}, considered above, deserves to and will be called a
representation by derivations of $A$ on $L$.
If $A \to M$ is a tangential Lie algebroid, then $A \ltimes L \to M$
will also be one.

\subsection{Differential operators}

We will from now on assume that $A$ denotes the vector bundle
determined by the structural Lie algebra $\VV$ and vice versa.

\begin{definition}
Let $\Diff{\VV}$ denote the algebra of differential operators
generated by $\VV$, where the vector fields are regarded as derivations
on functions.
\end{definition}

We also want to study differential operators with coefficients in
vector bundles.  Let $E_1 \to M$ and $E_2\to M$ be two vector
bundles. Embed $E_i \hookrightarrow M \times \CC^{N_i}$, $i=1,2$.
Denote by $e_i$ a projection in $M_{N_i}(\CI(M))$ whose (pointwise)
range is $E_i$. Then we define
\begin{equation}
        \Diff{\VV;E_1,E_2} := e_2
        M_{N_{2}\times N_{1}}(\Diff{\VV}) e_1.
\end{equation}
This definition of $ \Diff{\VV;E_1,E_2}$ is independent of the choices
of the embeddings $E_i \hookrightarrow M \times \CC^{N_i}$ and of the
choice of $e_i$. Elements of $\Diff{\VV;E_0,E_1}$ will be called {\em
differential operators generated by $\VV$}. In the special case $E_1 =
E_2= E$ we simply write $\Diff{\VV;E}$, the \emph{algebra of
differential operators on $E$ generated by $\VV$}.

It is possible to describe the differential operators in
$\Diff{\VV;E}$ locally on $M$ as follows.

\begin{lemma}\label{lemma.desc}\
A linear map $D : \Gamma(E) \to \Gamma(E)$ is in $\Diff{\VV;E}$ if,
and only if, for any trivialization $E \vert_U \cong U \times \CC^N$,
above some open subset $U \subset M$, the restriction $D\vert_U :
\Gamma(E\vert_U) \cong \CI(U) \otimes \CC^m \to \Gamma(E\vert_U)$ can
be written as a linear combination of compositions of operators of the
form $X \otimes 1$ and $f$, with $X \in \Gamma(A)$ and $f$ a smooth
endomorphism of the vector bundle $E \vert_U$.
\end{lemma}

\begin{proof}\
In a trivialization of $E$ above some open subset, we can assume
that $e$ is a constant matrix.
\end{proof}

\begin{example}\label{Example.1.18}
{\em De Rham differential generated by $\VV=\Ga(A)$.}\\ We define
for a section $\omega$ of $\Lambda^q A^*$
\begin{multline*}
        (d\omega)(X_0,\ldots, X_k) = \sum_{j=0}^q (-1)^j
        X_j(\omega(X_0,\ldots,\hat{X}_j, \ldots, X_k)) + \\
        \sum_{0 \le i < j \le q} (-1)^{i+j}
        \omega([X_i,X_j],X_0,\ldots,\hat{X}_i,
        \ldots,\hat{X}_j,\ldots, X_k).
\end{multline*}
This is well defined as $\Gamma(A)$ is closed under the Lie
bracket. By choosing a local basis of $A$ we see that this defines
a differential operator $\Gamma(\Lambda^q A^*) \to
\Gamma(\Lambda^{q+1}A^*)$ generated by $\VV=\Ga(A)$, the \emph{de
Rham differential}.

Assume now that $A\vert_{M_0} = TM_0$. The vector bundles
$\Lambda^q T^*M_0$ extend to bundles $\Lambda^q A^*$ on $M$. The
Cartan formula (e.g.\ \cite{BGV}) says that on $M_0$ the de Rham
differential is the de Rham differential of ordinary differential
geometry.
\end{example}

\begin{definition}\label{def.AvalC}\
Let $E \to M$ be a vector bundle. An {\em $A^*$-valued connection
on $E$} is a differential operator
\[
        D: \Gamma(E) \to \Gamma(E \otimes A^*)
\]
such that, for any $X \in \Gamma(A)$, the induced operator $D_X :
\Gamma(E) \to \Gamma(E)$ satisfies the usual properties of a
connection:
\begin{equation}
        (i)\ D_X(f\xi) = f D_X(\xi) + X(f) \xi; \quad \text{and}
        \quad (ii)\ D_{fX} \xi = f D_X(\xi).
\end{equation}
\end{definition}
It is clear from (i) that the operator $D_X$ is of first order.

Our definition of an $A^*$-valued connection is only slightly more
restrictive than that of $A$-connection introduced in
\cite{EvensLuW}. (In that paper, Evens, Lu, and Weinstein
considered (ii) only up to homotopy.)

Clearly if $D$ and $D'$ are $A^*$-valued connections on $E$ and,
respectively, $E'$, then $D'':= D \otimes 1 + 1 \otimes D'$ is an
$A^*$-valued connection on $E \otimes E'$.

See also \cite{eichorn3, eichorn4}

\section{Lie  structures at infinity\label{Sec.Unif.Str}}

In this section we introduce the class of manifolds with a Lie
structure at infinity, and we discuss some of their properties.
Our definition, Definition \ref{def.unif.str}, formalizes some
definitions from \cite{MelroseScattering}.

In some of the first papers on the analysis on open manifolds using
Lie algebras of vector, for example \cite{Cordes1, Cordes2,
CordesOwen1, TaylorG}, the vector fields considered were required to
vanish at infinity. In order to obtain more general results and in
agreement with the more recent papers on the subject (for example
\cite{CordesDoong, MelroseScattering, Muller, Vasy}), we do not make
this assumption. As a consequence, the comparison algebras that result
from our setting do not have in general the property that the
commutators are compact.

\subsection{Definition}

In the following, $\pa M$ denotes the union of all hyperfaces of a
manifold with corners $M$.

\begin{definition}\label{def.unif.str}\
A {\em Lie structure at infinity} on a smooth manifold $M_0$ is a
pair $(M,\VV)$, where
\begin{enumerate}[(i)]
\item $M$ is a compact manifold, possibly
with corners and $M_0$ is the interior of $M$;
\item $\VV$ is a structural Lie algebra of vector fields;
\item $\varrho_\VV: A_\VV\to TM$ induces an isomorphism on $M_0$, \ie\
$\varrho_\VV|_{M_0}:A_\VV|_{M_0}\to TM_0$ is a fiberwise
isomorphism.
\end{enumerate}
\end{definition}

If $M_0$ is compact without boundary, then it follows from the
above definition that $M = M_0$ and $A_\VV = TM$, so a Lie
structure at infinity on $M_0$ gives no additional information on
$M_0$. The interesting cases are thus the ones when $M_0$ is
non-compact.  Note that all the Examples \ref{bcal} -- \ref{theta}
are in fact Lie structures at infinity on the interior of $M$.

Here is now an explicit test for a Lie algebra of vector fields
$\VV$ on a compact manifold with corners $M$ to define a Lie
structure at infinity on the interior $M_0$ of $M$. This
characterization of Lie structures at infinity is in the spirit of
our discussion of local basis (see Equation \eqref{eq.local.basis}
and the discussion around it).

\begin{proposition}\label{prop.test}\ We have that the Lie algebra
$\VV \subset \Gamma(M; TM)$ defines a Lie structure at infinity on
$M_0$ if, and only if, the following conditions are satisfied:
\begin{enumerate}[(i)]
\item $\VV \subset \VV_b$, with $\VV_b$ defined in Equation
\eqref{eq.def.Vb}, and $\CI(M) \VV \subset \VV$.
\item If $x \in M_0$, $U$ is a compact neighborhood $U$ of $x$ in $M_0$,
and $Y$ is a vector field on $M_0$, then there exists $X \in \VV$,
such that $X|_U = Y|_U$.
\item \label{numthree}
If $x \in \pa M = M \smallsetminus M_0$, then we can find $n$
linearly independent vector fields $X_1, X_2, \ldots, X_n \in
\VV$, $n = \dim M$, defined on a neighborhood $U$ of $x$, such
that for any $X \in \VV$, there exist smooth functions $f_1,
\ldots, f_n \in \CI(U)$ uniquely determined by
\[
    X = \sum_{k = 1}^n f_k X_k \;\text{ on }\; U\,.
\]
\item There are functions $f_{ijk}\in \CI(U)$
(in particular smooth on the boundary $\pa M\cap U$) such that the
vector fields $X_j$ from \eref{numthree} satisfy $[X_i, X_j] =
\sum_{k = 1}^n f_{ijk} X_k$ on $U$.
\end{enumerate}
\end{proposition}

\begin{proof}\ The proof is an immediate translation of the
definition of a manifold with a Lie structure at infinity using the
description of projective $CI(M)$ module given at the end of
Subsection \ref{subs.pm} (especially Equation
\eqref{eq.local.basis}).
\end{proof}

\subsection{Riemannian manifolds with Lie  structures at infinity}

We now consider Riemannian metrics on $A \to M$.

\begin{definition}\label{def.R.usi}\
A manifold $M_0$ with a Lie structure at infinity $(M,\VV)$
together with a Riemannian metric on $A=A_\VV$, \ie\ a smooth
positive definite symmetric 2-tensor $g$ on $A$, is called a {\em
Riemannian manifold $M_0$ with a Lie structure at infinity}.
\end{definition}

In particular $g$ defines a Riemannian metric on $M_0$.  The
geometry of these metrics will be the topic of the next section.
Note that the metrics on $M_0$ that we obtain are not restrictions of
Riemannian metrics on $M$.
In the following section, we will prove for example
that $(M_0,g)$ is a complete Riemannian metric. Any curve joining a point on
the boundary $\partial M$ to the interior $M_0$ is necessarily
of infinite length.

\begin{example}
{\em Manifolds with cylindrical ends.}\\ A manifold $M$ with
cylindrical ends is obtained by attaching to a manifold $M_{1}$
with boundary $\pa M_1$ the cylinder $(-\infty,0] \times \pa M_1$,
using a tubular neighborhood of $\pa M_1$, where the metric is
assumed to be a product metric. The metric on the cylinder is also
assumed to be the product metric. Let $t$ be the coordinate of
$(-\infty,0]$.  By the change of variables $x = e^t$, we obtain
that $M$ is diffeomorphic to the interior of $M_1$ and $\VV =
\VV_b$. Other changes of variables lead us to different Lie
structures at infinity.  Similarly, products of manifolds with
cylindrical ends can be modeled by manifolds with corners and the
structural Lie algebra of vector fields $\VV_b$. This applies also
to manifolds that are only locally at infinity products of
manifolds with cylindrical ends.
\end{example}

\subsection{Bi-Lipschitz equivalence}
It turns out that the metric on a manifold with a Lie structure at
infinity is essentially unique, namely any two such metrics are
bi-Lipschitz equivalent (see the corollary below).

\begin{lemma}\label{lemma.bilip}
We assume that a manifold $M_0$ which is the interior of a compact
manifold with corners $M$ carries two Lie structures at infinity
$(M,\VV_1)$ and $(M,\VV_2)$ satisfying $\VV_1\subset \VV_2$.
Furthermore, let $g_j$ be Riemannian metrics on $A_{\VV_{j}}$,
$j=1,2$.  Then there is a constant $C$ such that
\[
        g_2(X,X)\leq C g_1(X,X)\qquad \mbox{for all}
        \quad X\in TM_0.
\]
\end{lemma}

\begin{proof}
The pull-back of $g_2$ to $A_{\VV_1}$ is a non-negative symmetric
two-tensor on $A_{\VV_1}$. The statement then follows from the
compactness of $M$.
\end{proof}

As a consequence the volume element of $g_2$ is bounded by a multiple
of the volume element of $g_1$. Furthermore, we have inclusions of
$L^p$-functions: $L^p(M_0,g_1)\hookrightarrow L^p(M_0,g_2)$.

\begin{corollary}\label{cor.bi.lip}
If two Riemannian metrics $g_1$, $g_2$ on $M_0$ are Riemannian metrics
for the same Lie structure at infinity $(M,\VV)$, then they are
bi-Lipschitz, \ie\ there is a constant $C>0$ with
\begin{equation*}
        C^{-1}g_2(X,X) \le g_1(X,X) \le C g_2(X,X)
\end{equation*}
for all $X\in TM_0$.
In particular, $C^{-1} d_2(x,y) \le d_1(x,y) \le C d_2(x,y)$, where
$d_i$ is the metric on $M_0$ associated to $g_i$.
\end{corollary}

\begin{proof}\
The first part is clear. The proof of the last statement is obtained
by comparing the metrics on a geodesic for one of the two metrics.
\end{proof}

\section{Geometry of Riemannian manifolds with Lie structures
at infinity\label{Sec.3}}

We now discuss some geometric properties of Riemannian manifolds with
a Lie structure at infinity. We begin with a simple observation about
volumes.

\subsection{Volume} Let
$d\vol$ be the volume element on $M_0$

\begin{proposition}\label{prop.vol}\ Let $M_0$ be a  Riemannian
manifold with Lie structure $(M,\VV,g)$ at infinity. Let $f \ge 0$ be
a smooth function on $M$. If $\int_{M_0}f d\vol < \infty$, then $f$
vanishes on each boundary hyperface of $M$.  In particular the volume
of any non-compact Riemannian manifold with a Lie structure at
infinity is infinite.
\end{proposition}

\begin{proof}\
Because of Lemma~\ref{lemma.bilip}, we can assume that $\Ga(A)$
are the vector fields tangential to the boundary.  For simplicity in
notation, let us assume that $M$ is a compact manifold with boundary.
Let $d\vol'$ be the volume element on $M$ associated to some metric on
$M$ that is smooth up to the boundary. Then $d\vol \ge C x^{-1}d\vol'$
for any boundary defining function $x$ and a constant $C$.

So, if $f$ is non-zero on $\pa M$ with defining function
$x$, then
$$
        \int_{M_0} f d\vol \ge \int_M f x^{-1}d\vol' = \infty.
$$
\end{proof}

\subsection{Connections and Curvature}
Most of the natural differential operators between bundles
functorially associated to the tangent bundle extend to differential
operators generated by $\VV$, with the tangent bundle replaced by
$A$. The main example is the Levi-Civita connection.

\begin{lemma}\label{lemma.LeviCivita}\ Let $M_0$ be a
Riemannian manifold with a Lie structure $(M,\VV,g)$ at infinity.
Then the Levi-Civita connection $\nabla : \Gamma(TM_0) \to \Gamma(TM_0
\otimes T^*M_0)$ extends to a differential operator in $\Diff{\VV;A,A
\otimes A^*}$, also denoted by $\nabla$. In particular, $\nabla_X \in
\Diff{\VV;A}$, for any $X \in \Gamma(A)$, and it satisfies
\begin{equation*}
        \nabla_X(fY) = X(f) Y + f \nabla_X(Y) \;\text{ and }\;
        X\<Y,Z\> = \< \nabla_XY,Z\> + \< Y,\nabla_XZ\>
\end{equation*}
for all $X, Y, Z \in \Gamma(A)$ and $f \in \CI(M)$. Moreover, the
above equations uniquely determine $\nabla$.
\end{lemma}

\begin{proof}\
Suppose $X,Y \in \VV = \Gamma(A) \subset \Gamma(TM)$.  We shall define
$\nabla_X Y$ on $M_0$ using the usual Levi-Civita connection $\nabla$
on $TM_0$. We need to prove that there exists $X_1 \in \Gamma(A)$
whose restriction to $M_0$ is $\nabla_X Y$.

Recall (for example from \cite{BGV}), that the formula for $\nabla_X
Y$ is given by
\begin{multline}
        2\<\nabla_X Y, Z\> = \<[X,Y],Z\> - \<[Y,Z],X\> + \<[Z,X],Y\>
        \\ + X\<Y,Z\> + Y\<Z,X\> - Z\<X,Y\>.
\end{multline}
Suppose $X,Y,Z \in \Gamma(A)$ in the above formula. We see then that
the function $2\<\nabla_X Y, Z\>$, which is defined a priori only on
$M_0$, extends to a smooth function on~$M$. Since the inner product
$\<\;,\;\>$ is the same on $A$ and on $TM_0$ (where they are both
defined), we see that the above equation determines $\nabla_X Y$ as a
smooth section of $A$. This completes the proof.
\end{proof}

The above lemma has interesting consequences about the geometry of
Riemannian manifolds with Lie structures at infinity.

Using the terminology of $A^*$-valued connections (see Definition
\ref{def.AvalC}), Lemma \ref{lemma.LeviCivita} can be formulated as
saying that the usual Levi-Civita connection on $M_0$ extends to an
$A^*$-valued connection on $A$. Similarly, we get $A^*$-valued
connections on $A^*$ and on all vector bundles obtained functorially
from $A$. We use this remark to obtain a canonical $A^*$-valued
connection on the bundles $A^{* \otimes k} \otimes \Lambda^2 A^*
\otimes \End(A)$.  (Here $E^{\otimes k}$ denotes $E \otimes \ldots
\otimes E$, $k$-times, as usual.)

We define the Riemannian curvature tensor as usual
\begin{equation}
        R(X,Y) := \nabla_X \nabla_Y - \nabla_Y \nabla_X -
        \nabla_{[X,Y]} \in \Gamma(\End(TM_0)),
\end{equation}
where $X,Y \in \Gamma(TM_0)$.  We will regard $R$ as a section of
$\Lambda^2 T^*M_0 \otimes \End(TM_0)$.  Then the covariant derivatives
$\nabla^k R \in \Gamma(T^*M_0 ^{\otimes k} \otimes \Lambda^2 T^*M_0
\otimes \End(TM_0))$ are defined.

\begin{corollary}
\label{cor1.20}\
If $M_0$ and $\nabla$ are as above (Lemma \ref{lemma.LeviCivita}),
then the Riemannian curvature tensor $R(X,Y)$ extends to an
endomorphism of $A$, for all $X,Y \in \Gamma(A)$. Moreover, each
covariant derivative $\nabla^k R$ extends to a section of $A^{*
\otimes k} \otimes \Lambda^2 A^* \otimes
\End(A)$ over $M$ and hence they are all bounded on $M_0$.
\end{corollary}

\begin{proof}\
Fix $X,Y \in \Gamma(A)$. Then $[X,Y] \in \Gamma(A)$ and hence
\begin{equation*}
        \nabla_X,\nabla_Y, \nabla_{[X,Y]} \in \Diff{\VV; A},
\end{equation*}
by Lemma \ref{lemma.LeviCivita}. It follows that $R(X,Y) \in
\Diff{\VV;A}$. Since $R(X,Y)$ induces an endomorphism of $TM_0$ and
$M_0$ is dense in $M$, it follows that $R(X,Y) \in \End(A)$.

Once we have obtained that $R \in \Gamma(\Lambda^2 A^* \otimes
\End(A))$, we can apply the $A^*$-valued Levi-Civita connection to
obtain
\begin{equation*}
        \nabla^k R \in \Gamma(A^{* \otimes k} \otimes \Lambda^2 A^*
        \otimes \End(A)).
\end{equation*}
The boundedness of $\nabla^k R$ follows from the fact that $M$ is
compact.
\end{proof}

The covariant derivative
\begin{equation}
        \nabla_X : \Gamma(A^{\otimes k} \otimes A^{* \otimes j}) \to
        \Gamma(A^{\otimes k} \otimes A^{* \otimes (j + 1)}), \quad X
        \in \Gamma(A),
\end{equation}
will be called, by abuse of notation, the {\em $A^*$-valued
Levi-Civita connection}, for any $k$ and $j$. Sometimes, when no
confusion can arise, we shall call this $A^*$-valued connection
$\nabla$ simply the {\em Levi-Civita connection}.

\subsection{Ehresmann connections}
\def\piN{\pi_{N}}
\begin{definition}\label{def.a.val}
Let $\piN:N\to M$ be a submersion of manifolds with corners, and let
$A \to M$ be a boundary tangential Lie algebroid.  Smooth sections of
the bundle $\bigwedge^p\thickpull{\piN}A\to N$ are called
\emph{$A^*$-valued $p$-forms on $N$}.
\end{definition}

The fiber of $A$ in $p\in M$ is denoted by $A_p$.

\begin{definition}\label{def.ehresmann}
An \emph{Ehresmann connection} on $\piN:N\to M$ with respect to $A$ is
a smooth field $x\mapsto\tau_x$, $x\in N$, of $n$-dimensional
subspaces of $\thickpull{\piN}A$ such that $(\piN)_* :
\thickpull{\piN} A \to A$ restricts to an isomorphism $\tau_x\to
A_{\piN x}$, for all $x\in N$, the inverse
$\left((\piN)_*\right)^{-1}: A_{\piN x}\to \tau_x$ is called a
\emph{horizontal lift}.
\end{definition}

We chose the terminology ``Ehresmann connection'' to honor important
work of Ehresmann's on the subject \cite{ehresmann}.

\begin{example}\
(a)\ Let $\pi:V\to M$ be a vector bundle. The definition of an
Ehresmann connection generalizes the notion of $A^*$-valued connection
(in the sense of vector bundle connections). In fact, let $\nabla$ be
an $A^*$-valued connection. Then we obtain the Ehresmann connection as
follows: For $X_0\in V$, $p=\pi(X_0)$ we extend $X_0$ to a local
section $X:U\to V$, where $U$ is a neighborhood of $p$ in $M$. We
define the horizontal lift
\begin{eqnarray*}
  H_{X_0}:A_p& \to & (\thickpull{\pi}A)_{X_0}\\ Y & \mapsto &
          (X_*)_p(Y)- \na_Y X,
\end{eqnarray*}
where $\na_Y X\in V_p\subset (\thickpull{\pi}A)_{X_0}$.  It is easy to
check that this map does not depend on the extension $X$ of $X_0$,
that $H_{X_0}$ is injective and that we have $(\pi_*)_{X_0}\circ
H_{X_0}=\id$.

Then $\tau_{X_0}:=\mathop{\mathrm{im}} H_{X_0}$ is an Ehresmann
connection on $V$.  The associated Ehresmann connection completely
characterizes the $A^*$-valued connection.  However, there are
Ehresmann connections on $V$ that do not come from $A^*$-valued
connections (they are not ``compatible'' with the vector space
structure).\vspace*{2mm}

\noindent (b)\ If the $A^*$-valued connection is metric with respect
to a chosen metric on $V$, then the Ehresmann connection is tangential
to the sphere bundle in $V$ with respect to that metric.
\end{example}

\subsection{Geodesic flow}
For any boundary tangential Lie algebroid $A$ equipped with a metric,
let $S(A)$ be the unit tangent sphere in $A$,
$$
        S(A):=\left\{ v \in A \,|\, \|v\|=1\right\}.
$$
The canonical projection map $\pi:S(A)\to M$ is a submersion of
manifolds with corners. Let $\thickpull{\pi}A$ be the thick pull-back
of $A$.

The manifold (with corners) $S(A)$ carries an Ehresmann connection and
a horizontal lift $H$ given by the Levi-Civita-connection on~$A$.

\begin{definition}
The \emph{geodesic spray} is defined to be the map
\begin{equation}
        S(A) \ni X \longmapsto H_X(X) \in \thickpull{\pi}A,
\end{equation}
which defines a section $S$ of $\thickpull{\pi}A\to S(A)$.  The flow
of this vector field is called the \emph{geodesic flow}.
\end{definition}

Restricted to the interior of the manifold, these concepts
recover the analogous concepts of ordinary Riemannian geometry.

By definition, the image of $S$ through the anchor map is a vector
field along $S(A)$ that is tangential to all the boundary faces of
$S(A)$. These boundary faces are preimages of the boundary faces of
$M$ under $\pi$.

\begin{lemma}
Let $A$ be a boundary tangential Lie algebroid and let $X\in
\Gamma(A)$. Then $X$ is complete in the sense that the flow lines
$\phi_t$ of $X$ are defined on $\RR$. The flow $\phi_t$ preserves the
boundary depth. In particular, flow lines emanating from $N_0:= N
\smallsetminus \pa N$ stay in $N_0$.
\end{lemma}

\begin{proof}\
For any boundary defining function $x_H$ one has
$$
    \left. \frac{d}{dt}\right|_{t=0}x_H(\phi_t)= dx_H(X)=0,
$$
hence the flow preserves the boundary depth.  In particular, the flow
preserves the boundary.  Let $I=(a,b)$ be a maximal open interval on
which one particular flow-line is defined. Let $t_i\to b$.  Assume
that $b < \infty$. Since $N$ is compact, after passing to a suitable
subsequence, we can assume that $\phi_{t_i}$ converges to $p\in N$.
In a neighborhood of $p$, the flow exists, which contradicts the
maximality of $b$. Hence $b=\infty$.  The proof for $a=-\infty$ is
similar.
\end{proof}

By applying this lemma to the geodesic flow on $N=S(A)$, we obtain two
corollaries.

\begin{corollary}\label{cor.complete}\
Let $M_0$ be a Riemannian manifold with a Lie structure $(M,\VV,g)$ at
infinity. Then $M_0$ is complete in the induced metric $g$.
\end{corollary}

\begin{corollary}
Let $M_0$ be a Riemannian manifold with a Lie structure $(M,\VV,g)$ at
infinity. Let $X\in A_p$, $p\in M$. Then the boundary depth of $\exp_p
A_p$ equals the boundary depth of $p$.
\end{corollary}

\subsection{Positive injectivity radius}
We already know that Riemannian manifolds with a Lie structure at
infinity are complete and have bounded sectional curvature.  For many
analytic statements it is very helpful if we also know that the
injectivity radius ${\rm inj}(M_0,g)= \inf_{p\in M_0} {\rm inj}_p$ is
positive.  For example see \cite{Kordyukov}, where a ``uniform bounded
calculus of pseudodifferential operators'' was defined on a manifold
with bounded geometry. Hebey \cite[Corollary~3.19]{hebey} proves
Sobolev embeddings for manifolds with bounded geometry, \ie\ complete
Riemannian manifolds with positive injectivity radius and bounded
covariant derivatives $\nabla^k R$ of the Riemannian curvature tensor
$R$.  We will say more about this in a future paper.

\begin{conjecture}\
All Riemannian manifolds with Lie structures at infinity have positive
injectivity radius.
\end{conjecture}

We now introduce two conditions on a Riemannian manifold $M_0$ with a
Lie structure at infinity $(M,\VV)$, see Definitions~\ref{def.lce}
and~\ref{def.cvfe}, and prove that if any of these conditions holds,
then the injectivity radius of $M_0$ is positive.

\begin{definition}\label{def.lce}\
A manifold $M_0$ with a Lie structure at infinity $(M,\VV)$ is said to
satisfy the \emph{local closed extension property for $1$-forms} if
any $p\in \partial M$ has a small neighborhood $U\subset M$ such that
any $\alpha_p\in A^*_p$ extends to a closed one-form on $U$.
\end{definition}

\begin{example}\
For the b-calculus, the local closed extension property holds, because
in the notation of Example~\ref{bcal} the locally defined closed
$1$-forms $dx^j/x^j$ and $dy^k$ span $A^*_p=(T^b_pM)^*$ for any $p\in
\partial M$.
\end{example}

\begin{theorem}\label{theo.lce.pir}\
Let $M_0$ be a manifold with a Lie structure at infinity $(M,\VV)$
which satisfies the local closed extension property.  Then for any
Riemannian metric $g$ on $A$ the injectivity radius of $(M_0,g)$ is
positive.
\end{theorem}

\begin{proof}\
We prove the theorem by contradiction.  If the injectivity radius is
zero, then, as the curvature is bounded, there is a sequence of
geodesics loops $c_i:[0,a_i]\to M_0$, parametrized by arc-length, with
$a_i\to 0$.  Because of the compactness of $S(A)$ we can choose a 
subsequence such that $\dot c_i(0)$ converges to a vector $v\in
S(A)$.  Obviously, the base point of $v$ has to be in $\partial M$.

By the local closed extension property, there is a closed 1-form
$\alpha$ on a small neighborhood of the base-point of $v$ such that
$\alpha(v)\neq 0$.  On the other hand, because of the closedness of
$\alpha$ we get for sufficiently large $i$
$$
         0 = \int_0^1 \alpha(\dot c_i(a_it))\, dt = \int_0^1
        \alpha(\phi_{a_i t}(\dot c_i(0))\, dt,
$$
where $\phi_t:S(A)\to S(A)$ denotes the geodesic flow.  As $i\to
\infty$, the integrand converges uniformly to $\alpha(v)$, thus we
obtain the contradiction $\alpha(v)=0$.
\end{proof}

In the remainder of this subsection we will prove another sufficient
criterion.

\begin{definition}\
Let $\phi:[0,\infty)^{n-k}\times \RR^k\to U\subset M$ be a local
parametrization of $M$, \ie $\phi^{-1}$ is a coordinate chart.  Then
for $v\in \RR^n$, the local vector field $\phi_*(v)$, \ie the image of
a constant vector field $v$ on $\RR^n$ is called a \emph{coordinate
vector field} with respect to $\phi$.
\end{definition}

\begin{definition}\label{def.cvfe}\
A manifold $M_0$ with a Lie structure at infinity $(M,\VV)$ is said to
satisfy the \emph{coordinate vector field extension property} if
$A_\VV$ carries a Riemannian metric $g$ such that for any $p\in
\partial M$ there is a parametrization $\phi: [0,\infty)^{n-k}\times
\RR^k\to U$ of a neighborhood $U$ of $p$ such that
\begin{enumerate}[(i)]
\item for any $v\in \RR^n\setminus\{0\}$ the normalized coordinate vector
field
$$
        X_v=\frac{\phi_*(v)}{\sqrt{g(\phi_*(v),\phi_*(v))}}
$$
which a priori is only defined on $U\cap M_0$
extends to a section of $A|_U$,
\item for linearly independent vectors $v$ and $w$, $X_v(p)$ and
$X_w(p)$ are linearly independent.\label{def.cvfe.two}
\end{enumerate}
\end{definition}
\noindent Note that Property (i) is equivalent to claiming that
$$
    f_v:=\frac{1}{\sqrt{g(\phi_*(v),\phi_*(v))}}
$$
extends to a smooth function on $M$.

\begin{theorem}\label{theo.cvfe.pir}\
Let $M_0$ be a manifold with a Lie structure at infinity $(M,\VV)$
that satisfies the coordinate vector field extension property.  Then
for any Riemannian metric $g$ on $A$, the injectivity radius of
$(M_0,g)$ is positive.
\end{theorem}

The theorem will follow right away from
Proposition~\ref{prop.crit.pir}, Lemma~\ref{lemma.cvfe.implies.con},
and Lemma~\ref{lem.con}, which we proceed to state and prove after the
following definition.

In the following, balls in euclidean $\RR^n$ will be called 
\emph{flat balls} or sometimes even just \emph{$\epsilon$-ball}.

\begin{definition}\ For $C\geq 1$ and $\ep>0$, we say that $(M_0,g)$ is
\emph{locally $C$-bi-Lipschitz to an $\ep$-ball} if each point $p\in
M_0$ has a neighborhood that is bi-Lipschitz diffeomorphic to a flat
ball of radius $\ep$ with bi-Lipschitz constant $C$.
\end{definition}

\begin{proposition}\label{prop.crit.pir}\
Let $(M_0,g)$ be a complete Riemannian manifold with
bounded sectional curvature. Then the following conditions are
equivalent:
\begin{enumerate}
\item $(M_0,g)$ has positive injectivity radius.
\item There are numbers $\de_1>0$ and $C>0$ such that any
loop of length $\de\leq \de_1$ is the boundary of a disk of
diameter $\leq C\cdot\de$.
\item There is $C>0$ and $\ep>0$ such that $(M_0,g)$ is locally
$C$-bi-Lipschitz to a ball of radius $\ep$.
\end{enumerate}
\end{proposition}

\begin{proof}\
(1)$\Longrightarrow$(3): Let $\rho>0$ be the injectivity radius of
$(M_0,g)$.  Then $B_{\rho/2}(p)$ is $C$-bi-Lipschitz to a flat $\epsilon$-ball
with $C$ independent from $p$.

(3)$\Longrightarrow$(2): Under the condition of (3) any loop of length
$\leq 2\ep/C$ based in $p$ lies completely inside $B_{\ep/C}(p)$. On
the other hand $B_{\ep/C}(p)$ is contained in a neighborhood $U$ of
$p$ which is $C$-bi-Lipschitz to a flat ball of radius $\ep$.  Inside
a flat ball any short loop is the boundary of a disk of small
diameter.  Hence (2) follows from (3) with $\de_1:=2\ep/C$ and with
the same $C$ in (2) as in (3).

(2)$\Longrightarrow$(1): Because sectional curvature is bounded from
above, there is a $\rho>0$ such that there are no conjugate points
along curves of length smaller than $\rho$.  For each $p\in M$, the
exponential map is a local diffeomorphism from $B_\rho(p)$ into~$M$.
We want to show that the exponential map is injective on any ball of
radius $\ep:=\min\{\de_1/2,\rho/(4C)\}$.  For this we assume that
$\exp_p(q_1)=\exp_p(q_2)$, $q_1,q_2\in B_\ep(0)\subset T_pM$.  Then
$$
        t\mapsto \begin{cases} \exp_p(tq_1)&\mbox{for } 0\leq t \leq
        1\\ \exp_p((2-t)\,q_2)&\mbox{for } 1\leq t \leq 2 \end{cases}
$$
is a closed loop of length $\leq 2\ep$.  Because of the conditions in
the proposition, this loop is the boundary of a disk of diameter $\leq
\rho/2$. Because the exponential map is a local diffeomorphism, it is
not difficult to see that such disks lift to $T_pM$. Hence $q_1=q_2$
which yields injectivity.

\komment{(1)$\Longrightarrow$(2): In this part of the proof we only
need a lower bound on sectional curvature.  Fix a number $\de_1$
smaller than the injectivity radius.  Let $c:[0,\de]\to M$ be a closed
loop parametrized by arc-length, $\de\leq \de_1$.  Then a
straightforward calculation in geodesic normal coordinates $\phi$
around $c(0)=c(\de)$ shows that the area of the image of
$$
        \Psi:[0,\de]\times [0,1]   \to M
$$
$$
        \Psi(t,s)= \phi^{-1}\left((1-s) \phi(c(0)) + s \phi(c(t))\right)
$$
is bounded from above by $\de^2+ C(\max |R|,\de_1) \de^3$ and the
diameter is bounded from above by $\de+ C(\max |R|,\de_1) \de^2$.
Hence the proposition holds.}
\end{proof}

As a consequence of this proposition, the property of having positive
injectivity radius is a bi-Lipschitz invariant inside the class of
complete Riemannian manifolds with bounded curvature.

\begin{corollary}\label{cor.bilipschitz.one}\
Let $g_1$ and $g_2$ be two complete metrics on $M_0$, such that there
is $C>0$ with
$$
        C^{-1} g_1 \le g_2 \le C g_1,
$$
$$
        |R_{g_1}|<C, \quad |R_{g_2}|<C.
$$
Then $(M_0,g_1)$ has positive injectivity radius if, and only if,
$(M_0,g_2)$ has positive injectivity radius.
\end{corollary}

\begin{proof}\
Proposition~\ref{prop.crit.pir} gives necessary and sufficient
criteria for positive injectivity radius, that are bi-Lipschitz
equivalent.
\end{proof}

Together with Corollary~\ref{cor.bi.lip} we obtain.

\begin{corollary}\label{cor.pos.inj.alg}\
Suppose that $M_0$ is a manifold with a Lie structure $(M,\VV)$ at
infinity and let $g$ and $h$ be two metrics on $A$.  Then
$(M_0,g|_{M_0})$ has positive injectivity radius if, and only if,
$(M_0,h|_{M_0})$ has positive injectivity radius.
\end{corollary}

\begin{definition}\label{def.controlled}\
Let $M_0$ be a Riemannian manifold with a Lie structure $(M,\VV,g)$ at
infinity. We say that the Lie structure at infinity is \emph{controlled} if for
all $p\in \partial M$ there is a parametrization
$\phi:[0,\infty)^{n-k}\times \RR^k\to U$ around $p$, a $\delta>0$ and
a constant $C>0$ such that for all $x\in M_0 \cap U$ and all $v\in
\RR^n$ the inequality
$$
        \max_{y\in B_\delta(x)}g_y\Bigl(\phi_*(v)_y,\phi_*(v)_y\Bigr)
        < C \min_{y\in
        B_\delta(x)}g_y\Bigl(\phi_*(v)_y,\phi_*(v)_y\Bigr)
$$
holds. Here $B_\de(x)$ denotes the ball of radius $\de$ around $x$
with respect to the metric~$g$.
\end{definition}

\begin{lemma}\label{lemma.cvfe.implies.con}\
Let $M_0$ be a manifold with a Lie structure at infinity $(M,\VV)$
satisfying the coordinate vector field extension property with the
metric $g$.  Then $(M,\VV,g)$ is controlled.
\end{lemma}

\begin{proof}\
Let $\phi$ be a parametrization of a neighborhood $U$ of $p\in \pa M$
and let $\phi_*(v)$ and $\phi_*(w)$ be arbitrary coordinate vector
fields.  Because of Property (i) in Definition~\ref{def.cvfe}, the
(local) functions $f_v:=1/\sqrt{g(\phi_*(v),\phi_*(v))}$ and
$f_w:=1/\sqrt{g(\phi_*(w),\phi_*(w))}$ extend to the boundary, and
$X_v:=f_v \phi_*(v)$ and $X_w:=f_w \phi_*(w)$ are sections of $A$.
For linearly independent $v$ and $w$ we calculate
\begin{eqnarray*}
        [X_v,X_w]& = & [f_v \phi_*(v),f_w \phi_*(w)]=
        X_v(f_w)\phi_*(w) -X_w(f_v)\phi_*(v)\\& = & X_v(\log f_w)X_w -
        X_w(\log f_v)X_v.
\end{eqnarray*}
This is again a section of $A$, and because of
Property~(\ref{def.cvfe.two}) in Definition~\ref{def.cvfe}, $X_v$ and
$X_w$ are even linearly independent on the boundary $\partial M$.  As
a consequence, $X_v(\log f_w)$ is bounded.  On the other hand, if
$v=w$, then because $f_w$ extends to the boundary $X_v(\log f_w)=f_v
\phi_*(v)(f_v)/f_v=\phi_*(v)(f_v)$ is also bounded. Hence, in both
cases:

\begin{equation}
|X_v(\log g(\phi_*(w),\phi_*(w)))| \leq C.\label{est.log}
\end{equation}

By summing up, one immediately sees that $C$ can be chosen
independently from the choice of $v$. Hence for fixed $w$,
$\eref{est.log}$ holds uniformly for any unit vector $X$.

We take two arbitrary points $y_1, y_2\in B_\de(x)$, $x\in U\cap M_0$.
We can assume that $B_\de(x)\subset U$.  We join $y_1$ and $y_2$ by a
path $c:[0,\tilde\de]\to M_0$, parametrized by arc-length,
$\tilde\de\leq 2 \de$. We estimate
\begin{eqnarray*}
  &&\left|\log\left(
  \frac{g_{c(0)}(\phi_*(w)_{c(0)},\phi_*(w)_{c(0)})} {g_{c(\tilde\de)}
  (\phi_*(w)_{c(\tilde\de)}, \phi_*(w)_{c(\tilde\de)})} \right)
  \right|\\ & = & \left|\int_0^{\tilde\de} \partial_{\dot c(t)}\,
  \left(\log
  g_{c(t)}(\phi_*(w)_{c(t)},\phi_*(w)_{c(t)})\right)\,dt\right|\\
  &\leq & \tilde\de\, C.
\end{eqnarray*}
Hence, the quotient is bounded by an expression that depends only on
$p,U,\de$ and global data.  Thus, $(M,\VV,g)$ is controlled.
\end{proof}

\begin{lemma}\label{lem.con}\
If the boundary tangential Lie algebroid is controlled, then there is
$C>0$ and $\ep>0$ such that $(M_0,g)$ is locally $C$-bi-Lipschitz to
an $\ep$-ball.
\end{lemma}

\begin{proof}\
On the ball $B_\delta(x)$ we regard two metrics: the original metric
$g$ and the metric $\tilde g$ with $g_x=\tilde g_x$ and $\tilde g_x$
is constant in the local coordinate chart.

These metrics are bi-Lipschitz on $B$ with a bi-Lipschitz constant
$C_1$.  Thus, the ball $\tilde B_x$ of radius $\de/C_1$ around $x$
with respect to the metric $\tilde g$ is a flat ball.  Hence, $(\tilde
B_x)_{x\in M_0}$ are neighborhoods that are uniformly bi-Lipschitz to
$\de/C_1$-balls.
\end{proof}

This completes the proof of Theorem \ref{theo.cvfe.pir}. We continue with
some examples and applications.

\begin{example}\
The examples~\ref{bcal} to \ref{theta} satisfy the coordinate vector
field extension property and hence have positive injectivity radius.
\end{example}

\begin{example}\
Let $M=[0,\infty)\times \RR\times \RR\ni(t,x,y)$. The vector fields
$X=t^2\partial_t$, $Y=e^{-C/t}\partial_x$ and $Z=e^{-C/t}\partial_y$
span a free $C^\infty(M)$-module that is closed under Lie-brackets.
Hence, it is a structural Lie algebra of vector fields $\VV$.  We
define a metric $g$ by claiming that $X,Y,Z$ are orthonormal, then
$(M,g,\VV)$ satisfies the coordinate vector field extension property.
\end{example}

\begin{example}\label{example.inf.rot}\
Let $M=[0,\infty)\times \RR\times \RR\ni(t,x,y)$. The vector fields
$X=t^2\partial_t$, $Y=e^{-C/t}\big(\sin(1/t)\partial_x+
\cos(1/t)\partial_y\big)$ and $Z=e^{-C/t}\big(\cos(1/t)\partial_x-\sin
(1/t)\partial_y\big)$ span a free $C^\infty(M)$-module that is closed
under Lie-brackets.  Hence it is a structural Lie algebra of vector
fields $\VV$.  We define a metric by claiming that $X,Y,Z$ are
orthonormal.  In this example the normalized coordinate vector fields
$e^{-C/t}\partial_x$ and $e^{-C/t}\partial_y$ are {\em not} contained
in $\VV$.  Hence $(M,\VV)$ does not satisfy the coordinate vector
field extension property.  However, $(M_0:=M\backslash \partial M,g)$
is isometric to the interior of the previous example.
Corollary~\ref{cor.bilipschitz.one} together with
Theorem~\ref{theo.cvfe.pir} will show that $(M_0,g)$ has positive
injectivity radius, although the conditions in
Theorem~\ref{theo.cvfe.pir} are not satisfied for $(M,\VV)$ directly.
\end{example}

\subsection{Adjoints of differential operators}\
We shall fix in what follows a metric on $A$, thus we obtain a
Riemannian manifold $(M_0,g)$ with a Lie structure at infinity, which
will remain fixed throughout this section.

Let us now discuss adjoints of operators in $\Diff{\VV}$. The metric
on $M_0$ defines a natural volume element $\mu$ on $M_0$, and hence it
defines also a Hilbert space $L^2(M_0, d\mu)$ with inner product
$(g_1,g_2):= \int_{M_0} g_1 \overline{g_2} d\mu$. The formal adjoint
$\sh{D}$ of a differential operator $D$ is then defined by the formula
\begin{equation}
        (D g_1, g_2) = (g_1, \sh{D}g_2), \quad \forall g_1,g_2 \in
        \CIc(M_0).
\end{equation}
We would like to prove that $ \sh{D} \in \Diff{\VV}$ provided that $D
\in \Diff{\VV}$. To check this, we first need a Lemma. Fix a local
orthonormal basis $X_1, \ldots, X_n$ of $A$ (on some open subset of
$M$). Then $\nabla_{X_i}X = \sum c_{ij}(X) X_j$, for some smooth
functions $c_{ij}(X)$. Then $\operatorname{div}(X) := -\sum c_{jj}(X)$
is well defined and gives rise to a smooth function on~$M$. See
\cite{gallot.hulin.lafontaine}, Chapter IV. A.

\begin{lemma}\
Let $X \in \Gamma(A)$ and $f \in \CIc(M_0)$. Then
\begin{equation*}
        \int_{M_0} X(f) \mu =  \int_{M_0} f \operatorname{div}(X)
        \mu.
\end{equation*}
In particular, the formal adjoint of $X$ is $\sh{X} = - X +
\operatorname{div}(X) \in \Diff{\VV}$.
\end{lemma}

\begin{proof}
We know (\cite{gallot.hulin.lafontaine}, Example 4.6) that
$$
        \operatorname{div}(fX)= f \operatorname{div}(X) - X(f).
$$
The divergence theorem (e.g. \cite{gallot.hulin.lafontaine}, Chapter
IV.A.)  states for $X\in \Gamma(A)$ and compactly supported functions
$f$
\begin{equation*}
        0 = \int_{M_0} \operatorname{div}(fX)\,\mu = \int_{M_0} f
        \operatorname{div} (X) \,\mu - \int_{M_0} X(f)\,\mu ,
\end{equation*}
so
\begin{equation*}
        \int_{M_0} X(f) \mu = \int_{M_0} f \operatorname{div}(X) \mu.
\end{equation*}
Now, if we set $f=g_1 \overline{g_2}$, we see directly that
\begin{equation*}
        \int_{M_0} g_1 \overline{X(g_2)}\,\mu + \int_{M_0} X(g_1)
        \overline{g_2}\, \mu = \int_{M_0} X(g_1 \overline{g_2})\,\mu =
        \int_{M_0} \operatorname{div}(X) g_1 \overline{g_2} \,\mu.
\end{equation*}
This implies the formula for the adjoint of $X$.
\end{proof}

\begin{corollary}\
Let $(M_0,g)$ be a Riemannian manifold with a Lie structure at
infinity. The algebra $\Diff{\VV}$ is closed under taking formal
adjoints. Similarly, if $E$ is a hermitian vector bundle on $M$, then
$\Diff{\VV;E}$ is also closed under taking formal adjoints.
\end{corollary}

\begin{proof}\ The formal adjoint of a vector field $X \in \VV$,
when regarded as a differential operator on $M_0$, is given by $\sh{X}
= - X + \operatorname{div}(X)$. The adjoint of $f \in \CI(M)$ is given
by $\sh{f} = \overline{f}$.  Since $\Diff{\VV}$ is generated as an
algebra by operators of the form $X$ and $f$, with $X$ and $f$ as
above, and $\sh{(D_1D_2)} = \sh{D_2} \sh{D_1}$, this proves that
$\Diff{\VV}$ is closed under taking adjoints.

If $E$ is a hermitian vector bundle, then we can choose the embedding
$E \to M \times \CC^N$ to preserve the metric. Then the projection $e$
onto the range of $E$ is a self-adjoint projection in
$M_N(\Diff{\VV})$. The equation $e^*=e$ satisfied by $e$ guarantees
that $\Diff{\VV ; E} := e M_N(\Diff{\VV}) e$ is also closed under
taking formal adjoints.
\end{proof}

Similarly, we obtain the following easy consequence.

\begin{corollary}\label{cor.last}\
If $E_0,E_1 \to M$ are hermitian vector bundles, then the adjoint of
an operator $P \in \Diff{\VV;E_0,E_1}$ is in $\Diff{\VV;E_1,E_0}$.
\end{corollary}

\begin{proof}\
Write $E := E_0 \oplus E_1$ and use the resulting natural matrix
notation for operators in $\Diff{\VV;E}$.
\end{proof}

\section{Geometric operators\label{Sec.Geom.Op}}

In this section we will see that the Hodge Laplacian $(d+d^*)^2$ on
forms and the classical Dirac operator on a Riemannian (spin) manifold
$M_0$ with a Lie structure at infinity $(M,\VV)$ are differential
operators generated by $\VV=\Ga(A)$.  (See also \cite{LN1} for some
similar results).

Both the classical Dirac operator and $d+d^*$ are generalized Dirac
operators.  We will show that any generalized Dirac operator is a
differential operator generated by $\VV$.  Our approach follows
closely that in \cite{gromov.lawson}.

\subsection{Hodge-Laplacians} Recall from Example \ref{Example.1.18} that the
de Rham differential defines an element
$d \in \Diff{\VV; \Lambda^p A^* , \Lambda^{p+1} A^*}$.

\begin{proposition}\ On a Riemannian manifold $M_0$ with a Lie structure
at infinity $(M,\VV)$, the Hodge-Laplace operator
\begin{equation}
        \Delta_p = d^* d + d d^*=(d+d^*)^2\in \Diff{\VV; \Lambda^p A^*}
\end{equation}
that is, it is a differential operator generated by $\VV$.
\end{proposition}

\begin{proof}\
This follows directly from Corollary \ref{cor.last} and the
construction in Example \ref{Example.1.18}.
\end{proof}

\subsection{Principal bundles and connection-$1$-forms}
\def\fg{\mfk g} \def\pip{\pi_P} Let $E\to M$ be a vector bundle of
rank $k$ carrying a metric and an orientation.  In this subsection, we
will show that giving a metric $A^*$-valued connection on $E$ is
equivalent to giving a $A^*$-valued connection-$1$-form on the frame
bundle of $E$.  Our approach generalizes the case of Riemannian
manifolds (see e.g. \cite[II, \S 4]{LM}), hence we will omit some
details.

For simplicity, we will assume from now on that the vector bundle $A
\to M$ is orientable.

Let $P$ be a principal $G$-bundle. We denote the Lie algebra of $G$
with $\fg$.  The most important example will be the bundle of oriented
orthonormal frames of the bundle $E$, denoted by $\pip:\PSO(E)\to M$,
which is a principal $\SO(k)$-bundle.  Differentiating the action of
$G$ gives rise to the canonical map
$$\fg\to \Ga(TP), \quad V\mapsto \widetilde V.$$

\begin{definition}
An $A^*$-valued connection-$1$-form $\om$ is an $\fg\otimes
A^*$-valued $1$-form on $\PSO(A)$ satisfying the compatibility
conditions
$$
        \om(\widetilde{V})=V \quad \text{and} \quad g^*\om=\Ad_{g^{-1}} \om
        \quad \mbox{for all }\; V\in \so(\fg).
$$
If $\fg\subset \so(k)$, we write $\om=(\om_{ij})$ with respect to the
standard basis of $\so(k)$.  In particular, the $\om_{ij}$ are
$A^*$-valued $1$-forms on $\PSO(A)$ satisfying $\om_{ij}=-\om_{ji}$.
\end{definition}

Here ``$A^*$-valued'' is in the sense of Definition \ref{def.a.val}.
Any $A^*$-valued con\-nection-$1$-form on $P$ gives rise to a
$G$-invariant Ehresmann connection on the bundle $P$ via $\tau=\{V\in
\thickpull{\pip}A\,|\,\om(V)=0\}$.  It is easy to check that this
yields a one-to-one correspondence between $G$-invariant Ehresmann
connections and connection-$1$-forms.

\begin{proposition}\label{prop.con.eq}
Let $E\to M$ be a vector bundle. For any $A^*$-valued
connection-$1$-form on $\PSO(E)$ there is a unique metric connection
on $E$ satisfying the formula
$$
        \nabla e_i= \sum_{j=1}^n {\cal E}^*{\om}_{ji}\otimes e_j
$$
where ${\cal E}=(e_1,\dots,e_n)$ is a local section of $\PSO(E)$.
Conversely, any metric connection on $E$ arises from such an
$A^*$-valued connection-$1$-form.
\end{proposition}

Note that ${\cal E}^*{\om}_{ij}$ is a well-defined $A^*$-valued
$1$-form on $M$.

The proof is straightforward and runs completely analogous to
\cite[II, Proposition~4.4]{LM} with ordinary $1$-forms replaced by
$A^*$-valued $1$-forms.  As a result, we conclude that the Levi-Civita
connection on $A$ determines an $\SO(n)$-invariant Ehresmann
connection and an $A^*$-valued connection-$1$-form on $\PSO(A)$.

\subsection{Spin structures and spinors}
The results of the previous subsection now allow us to define the
classical Dirac operator in a coordinate-free definition manner.
\def\Pspin{P_{\rm Spin}}

\begin{definition}
A \emph{spin structure} on $(M,\VV)$ is given by a $Spin(n)$-principal
bundle $\Pspin$ over $M$ together with bundle map $\theta:\Pspin\to
\PSO(A)$ that is $\Spin(n)\to\SO(n)$-equivariant.
\end{definition}

The (thick) pull-back of any $\SO(n)$-invariant Ehresmann connection
on the principal $\SO(n)$ bundle $\PSO(A)\to M$ with respect to $A$
defines a $\Spin(n)$-invariant Ehresmann connection on $\Pspin\to M$
with respect to $A$. Similarly, by using the standard identification
of the Lie algebra of $\SO(n)$ with the Lie algebra of $\Spin(n)$, any
$A^*$-valued connection-$1$-form on $\PSO(A)$ pulls back to an
$A^*$-valued connection-$1$-form on $\Pspin$.

\begin{definition}\label{def.ext.p.bundle}
Let $P$ be a principal bundle with respect to the Lie group $G$. Let
$\fg$ and $\mfk h$ denote the Lie algebras of $G$ and $H$.  Let
$\rho:G\to H$ be an inclusion of Lie groups. Then any $A^*$-valued
connection-$1$-form on $P$ defines an \emph{induced} $A^*$-valued
connection-$1$-form on $P\times_\rho H$ via the formula
$$
        [X;Y] \mapsto \rho_*(\om(X))+Y \qquad\mbox{for all $X\in TP$
        and $Y\in TH$}
$$
\end{definition}
This definition does not depend on the choice of representative $[X;Y]$ as
the map is invariant under the action of $G$ on $P\times H$.

Now, let $\si_n:\Spin(n)\to \SU(\Si_n)$ be the complex spinor
representation, e.g. the restriction of an odd irreducible complex
representation of the Clifford algebra on $n$-dimensional space
\cite{LM}. The complex dimension of $\Si_n$ is $d_n:=2^{[n/2]}$.

\begin{definition}\label{def.spinor.bundle}\
Let $M_0$ be an $n$-dimensional Riemannian manifold with a Lie
structure at infinity $(M,\VV,g)$ carrying a spin structure
$\Pspin(A)\to\PSO(A)$.  The \emph{spinor bundle} is the associated
vector bundle $\Si M:=\Pspin(A)\times_{\si_n} \Si_n$ on $M$.
\end{definition}

Any metric $A^*$-valued connection on $A$ gives rise to an
$A^*$-valued connection on $\Si M$ as follows:
Proposition~\ref{prop.con.eq} defines an $A^*$-valued
connection-$1$-form on $\PSO(A)$ which can be pull-backed to
$\Pspin(A)$. With Definition~\ref{def.ext.p.bundle} applied to
$\rho=\si_n:\Spin(n)\to\SU(d_n)\subset \SO(2d_n)$ we obtain an
$A^*$-valued connection-$1$-form on $\Pspin(A)\times_{\si_n}\SO(2d_n)$
compatible with complex multiplication. Another application of
Proposition~\ref{prop.con.eq} yields a complex $A^*$-valued connection
on $\Si M$.

In particular, the Levi-Civita-connection on $A$ defines then a metric
connection on $\Si M$, the so-called \emph{Levi-Civita-connection}.

Recall that the spinor representation $\Si_n$ admits a
$\Spin(n)$-equivariant linear map
$$
        \RR^n \otimes \Si_n \to \Si_n: \quad X
        \otimes \phi \mapsto X \cdot \phi.
$$
satisfying
\begin{equation}\label{cliff.rel}
        \big( X \cdot Y  + Y \cdot X  + 2 g(X,Y) \big) \cdot \phi = 0
\end{equation}
for all $X,Y \in \RR^n$ and all $\phi \in \Si_n$, the so-called
\emph{Clifford multiplication relations}.  By forming the associated
bundles this gives rise to a bundle map $A\otimes \Si M\to \Si M$,
called \emph{Clifford multiplication}. Equation (\ref{cliff.rel}) is
satisfied for all $X, Y \in A, \phi \in \Si M$ in the same base point.

\subsection{Generalized Dirac operators}\label{subsec.gen.dirac}\

\noindent We now discuss Clifford modules in our setting.

\begin{definition}\label{def.Clifford.mod}\
A \emph{Clifford module} over $M$ is a complex vector bundle $W\to M$
equipped with a positive definite product $\<\cdot, \cdot\>$,
anti-linear in the second argument, an $A^*$-valued connection
$\nabla^W$, and a linear bundle map $A\otimes W\to W$, $X\otimes
\phi\mapsto X\cdot \phi$ called \emph{Clifford multiplication} such
that
\begin{enumerate}[(1)]
\item \label{def.Clifford.mod.one}
$$
        \big( X \cdot Y + Y \cdot X + 2 g(X,Y) \big) \cdot \phi = 0
$$
\item \label{def.Clifford.mod.two}
$\nabla^W$ is metric, \ie\
$$
        \partial_X\<\psi, \phi\>= \<\nabla^W_X \psi, \phi\> +
        \<\psi,\nabla^W_X\phi\>
$$
for $X\in \Gamma(A)$, $\phi,\psi\in \Gamma(W)$,
\item \label{def.Clifford.mod.three} Clifford multiplication with
vectors is skew-symmetric, \ie\
$$
        \<X\cdot \psi, \phi\> =\<\psi, X\cdot\phi\>
$$
for
$\phi,\psi\in \Gamma(W)$, $X\in \Gamma(A)$,
\item \label{def.Clifford.mod.four}
Clifford multiplication is parallel, \ie\
$$
        \nabla^W_X(Y\cdot \phi)= (\nabla^W_X Y)\cdot \phi
        +Y \cdot (\nabla^W_X \phi)
$$
for  $X\in \Gamma(A)$, $Y\in \Gamma(A)$, $\phi\in \Gamma(W)$.
\end{enumerate}
\end{definition}

The \emph{generalized Dirac operator} associated to a Clifford module $W$
is the first order operator $\Dir^W$ obtained by the following composition:
$$
        \Gamma(W) \stackrel{\nabla^W}{\longrightarrow}
        \Gamma(W \otimes A^*) \stackrel{\id \otimes \#}
        {\longrightarrow} \Gamma(W\otimes A)
        \stackrel{\cdot}{\longrightarrow} \Gamma(W)
$$
$$  \Dir^W:= \cdot \circ (\id \otimes \#) \circ \nabla^W.
$$
The last map is Clifford multiplication and $\#:A^*\to A$ is the
isomorphism given by $g$.

The principal symbol of any generalized Dirac operator is elliptic, as
for any non zero vector $X$, Clifford multiplication $X\cdot$ is an
invertible element of $\End(\Si M)$.

\begin{example}\
For any $p\in M$, we define the Clifford algebra $\Cl(A_p)$ as the
universal commutative algebra generated by $A_p$ subject to the
relation
$$
        X \cdot Y +Y \cdot X + 2 g(X,Y)1 = 0.
$$
Let $\Cl(A)$ be the Clifford-bundle of $(A,g)$, \ie the bundle whose
fibers at the point $p\in M$ is the Clifford algebra $\Cl(A_p)$.  The
$A^*$-valued connection on $A$ extends to an $A^*$-valued connection
on $\Cl(A)$.  Let $W=\Cl(A)$, equipped with the module structure given
by left multiplication. After identifying with the canonical
isomorphism $\Cl(A)\cong\Lambda^*(A)$, $e_{i_1}\cdot \ldots\cdot
e_{i_k}\mapsto e^b_{i_1}\wedge \ldots\wedge e^b_{i_1}$ for an
orthonormal basis $(e_i)$ with dual $(e^b_i)$, the generalized Dirac
operator on this bundle is the de Rham operator $d + d^*$.
\end{example}

\begin{example}\
If $M$ is spin then the spinor bundle from Definition~\ref{def.spinor.bundle}
is also a Clifford module.  The corresponding Dirac
operator is called the \emph{(classical) Dirac operator}.
\end{example}

\begin{example}
If $M$ is even dimensional and if $A$ carries a K\"ahler structure, 
then the Dolbeault operator 
$\sqrt{2}\left(\overline{\partial} +\overline{\partial}^*\right)$ acting on
$(0,*)$-forms is a generalized Dirac operator.
\end{example}

\begin{example}
Let $W$ be any complex vector bundle with a positive definite scalar
product and a Clifford multiplication such that
(\ref{def.Clifford.mod.one}) and (\ref{def.Clifford.mod.three}) of
Definition~\ref{def.Clifford.mod} are satisfied. Then there is always
a connection $\na^W$ on $W$ satisfying the compatibility conditions
(\ref{def.Clifford.mod.two}) and (\ref{def.Clifford.mod.four}). This
can be seen as follows: If $M$ is spin then we can write
$$
    W=\Si M \otimes V
$$
where $V$ is isomorphic to the homomorphisms from $\Si M$ to $W$ that
are $\Cl(A)$ equivariant. $V$ carries a compatible metric. After
choosing any metric $A^*$-valued connection $\na^V$ on $V$, the
product connection on $W$ satisfies (\ref{def.Clifford.mod.two}) and
(\ref{def.Clifford.mod.four}).

If $M$ is not spin, the connection can be constructed locally on a
open covering in the same way, and the connection can then be glued
together by using a partition of unity, hence we obtain the statement.
\end{example}

For any two sections $\si_1$ and $\si_2$ of $W$, we let
$$
        (\si_1,\si_2):=\int_M \< \si_1,\si_2.\>
$$
This expression is not always defined. However, it is well-defined
scalar product on generalized $L^2$-spinor fields, \ie\ generalized
spinor fields with $\int_M \<\si_i,\si_i\><\infty$.  It is also
well-defined, if one of the sections $s_i$ or $s_j$ has compact
support and the other is locally $L^2$.

For the benefit of the reader, let us recall the following basic
result (see for example \cite{gromov.lawson}).

\begin{proposition}[\cite{gromov.lawson}]
Generalized Dirac operators $\Dir$ on complete Riemannian manifolds
are formally self-adjoint and essentially self-adjoint.  More
concretely, for smooth sections $\si_i$ we have
$$
        (\Dir\si_1,\si_2)=(\si_1,\Dir\si_2)
$$
if at least one of the sections $\si_1$ or $\si_2$ has compact
support, and the maximal and minimal extension of $\Dir$ coincide,
hence $\Dir$ extends uniquely to a self-adjoint operator densely
defined on the $L^2$-sections of $W$.
\end{proposition}

For any choice of a connection as in the above theorem, the resulting
Dirac operator is generated by $\VV$.

\begin{theorem}\label{theorem.Dirac}\
Let $W \to M$ be a Clifford module.  Then the Dirac operator on $W$ is
generated by $\VV$:
$$
        \Dir^W \in \Diff{\VV;W}.
$$
\end{theorem}

\begin{proof}
The Dirac operator is the composition of Clifford multiplication and
the $A^*$-valued connection $\na^W$ on $W$.  Clifford multiplication
is a zero order differential operator generated by $\VV$.  The
$A^*$-valued connection $\na^W$ on $W$ is a first order differential
operator generated by $\VV$.  Hence the Dirac operator is also a first
order differential operator generated by $\VV$.
\end{proof}


\end{document}